\documentclass[12pt]{article}
\usepackage{amsmath,amssymb,amsthm,enumerate,verbatim}
\usepackage{eucal}

\newtheorem{thm}{Theorem}[section]
\newtheorem*{thm*}{Theorem}
\newtheorem{lem}[thm]{Lemma}
\newtheorem{prop}[thm]{Proposition}
\newtheorem*{prop*}{Proposition}
\newtheorem{cor}[thm]{Corollary}
\theoremstyle{remark}
\newtheorem*{rmk*}{Remark}
\newtheorem{rem}{Remark}[section]
\newtheorem{ex}[rem]{Example}
\theoremstyle{definition}
\newtheorem*{dfn*}{Definition}
\numberwithin{equation}{section}

\newtheorem*{namedprop}{\theoremname}
\newcommand{\theoremname}{testing}
\newenvironment{prop1}[1]{\renewcommand{\theoremname}{#1}
    \begin{namedprop}}
    {\end{namedprop}}

\newcommand{\norm}[1]{\left\Vert#1\right\Vert}

\newcommand{\abs}[1]{\left\vert#1\right\vert}
\newcommand{\set}[1]{\left\{#1\right\}}
\newcommand{\brac}[1]{\left(#1\right)}
\newcommand{\scalar}[1]{\left \langle #1 \right \rangle}

\newcommand{\Real}{\mathbb{R}}

\newcommand{\eps}{\varepsilon}

\def \RR {\mathbb R}

\def \eps {\varepsilon}

\def\mes{\text{mes}}

\begin{document}

\title{An isoperimetric inequality for uniformly log-concave
measures and uniformly convex bodies}

\author{
Emanuel Milman\textsuperscript{1} and Sasha
Sodin\textsuperscript{2}} \footnotetext[1]{Mathematics Department,
The Weizmann Institute of Science, Rehovot 76100, Israel. Email:
emanuel.milman@gmail.com. Supported in part by ISF and BSF.}
\footnotetext[2]{School of Mathematics, Raymond and Beverly Sackler
Faculty of Exact Sciences, Tel Aviv University, Tel Aviv, 69978,
Israel. Email: sodinale@tau.ac.il. Research was supported in part by
the European Network PHD, MCRN--511953.}

\maketitle

\begin{abstract}
We prove an isoperimetric inequality for the uniform measure on a
uniformly convex body and for a class of uniformly log-concave
measures (that we introduce). These inequalities imply (up to
universal constants) the log-Sobolev inequalities proved by
Bobkov--Ledoux \cite{BobkovLedoux} as well as the isoperimetric inequalities
due to Bakry-Ledoux \cite{BakryLedoux} and Bobkov--Zegarlinski 
\cite{BobkovZegarlinski}. We also recover a concentration inequality
for uniformly convex bodies, similar to that proved by
Gromov--Milman \cite{Gromov-Milman}.
\end{abstract}

\section{Introduction}

Let $V = (\RR^n, \| \cdot \|)$ be a normed space, and let $\mu$ be a
probability measure on $V$ with density $f = \exp(-g)$, $g: \RR^n
\to \RR \cup \{+ \infty\}$. If $g$ is convex, the function $f$ and
the measure $\mu$ are called log-concave. Log-concave functions and
measures boast many important properties (cf.\ Borell
\cite{Borell-logconcave}, Bobkov \cite{BobkovGaussianIsoLogSobEquivalent} et cet.)

In this note, we study more restricted classes of measures. Let
\[ \delta: \RR_+ \to \RR_+ \cup \{ + \infty \}, \]
and consider the following condition:
\begin{equation}\label{ulc}
\frac{g(x) + g(y)}{2} - g\brac{\frac{x+y}{2}} \geq \delta(\norm{x -
y})~.
\end{equation}

\begin{ex}
The log-concavity condition corresponds to $\delta \equiv 0$.
\end{ex}

By analogy with uniformly convex bodies (cf.\
Subsection~\ref{sec:intro-uniformly-convex-bodies}), we define the
modulus of convexity $\delta_{g,\norm{\cdot}}$ of $g$ with respect
to the norm $\norm{\cdot}$ as:
\[
\delta_{g,\norm{\cdot}}(t) := \inf \set{\frac{g(x) + g(y)}{2} -
g\brac{\frac{x+y}{2}} ; \norm{x-y} \geq t \text{ and } g(x),g(y) <
\infty}.
\]
If $\delta_{g,\norm{\cdot}}(t) > 0$ for all $t>0$, we say that $f$
and $\mu$ are uniformly log-concave, and that $g$ is uniformly
convex. Obviously, this notion does not depend on the choice of the
norm $\norm{\cdot}$.

It is easy to check that $\delta_{g,\norm{\cdot}}(t)/t$ is always a non-decreasing
function of $t$; therefore in the sequel we consider measures $\mu$ satisfying (\ref{ulc})
with respect to a function $\delta$ such that
\begin{equation}\begin{cases}\label{eq:delta-assumption}
\delta(t) > 0, &t > 0 \\
t \mapsto \delta(t)/t &\text{is non-decreasing.}
\end{cases} \end{equation}

\begin{ex}\label{e2}
Let $\| \cdot \| = | \cdot |$ be the Euclidean norm, and let
$\delta(t) = t^2/8$. Then (\ref{ulc}) holds iff $\mu$ has
log-concave density with respect to the standard Gaussian measure;
in other words, if $\mu$ satisfies the Bakry--\'Emery
curvature-dimension condition $CD(1, +\infty)$ (cf.\ Bakry and
\'Emery \cite{BakryEmery}; recall that the usual log-concavity of
$\mu$ is equivalent to $CD(0, +\infty)$).
\end{ex}

\begin{rem}\label{rem:aff}
The condition (\ref{ulc}) is translation invariant. Therefore one may extend it to
measures on an affine space $\mathbb{A}^n$ on which $V$ acts by translations; note
that both sides of (\ref{ulc}) are still defined. This point of view will be convenient
in Section~\ref{sec:isop}.
\end{rem}

\setcounter{subsection}{-1}
\subsection{Assumptions and Notation}

Unless mentioned otherwise, the sets in this note are Borel subsets
of $\RR^n$, and the measures are Borel measures on $\RR^n$.

The Lipschitz norm of a map $T: V_1 \rightarrow V_2$ between two normed
spaces $V_i = (X_i,\norm{\cdot}_i)$, $i=1,2$, is defined as:
\begin{equation} \label{eq:Lip-def}
\norm{T}_\text{Lip} = \sup_{x,y \in X_1, x\neq y}
\frac{\norm{T(x)-T(y)}_2}{\norm{x-y}_1}.
\end{equation}
$T$ is called Lipschitz if $\norm{T}_{Lip} < \infty$.
If
\[  \sup_{x,y \in K, x\neq y} \frac{\norm{T(x)-T(y)}_2}{\norm{x-y}_1} < +\infty \]
for any compact subset $K \subset X_1$, $T$ is called locally
Lipschitz.

A Borel map $T: V_1 \rightarrow V_2$ is said to push a measure $\mu$
on $V_1$ forward to a measure $\lambda$ on $V_2$ (notation:
$T_*\mu = \lambda$) if $\mu(T^{-1}(B)) = \lambda(B)$ for every
$B \subset X_2$.

If $\mu$ is a probability measure on $V = (X,\norm{\cdot})$, the
Minkowski boundary measure associated with $\mu$ (and $\| \cdot \|$)
is defined by:
\begin{equation}
\mu^+_{\| \cdot \|} (A)
    = \liminf_{\eps \to 0} \frac{\mu(A_{\eps, \| \cdot \|})
        - \mu(A)}{\eps}~, \quad A \subset X,
\end{equation}
where
\[
A_{\eps,\| \cdot \|}
    = \left\{ x \in X \, \mid \, \exists y \in A, \| x - y \| < \eps \right\}
\]
is the $\eps$-extension of $A$ in the metric induced by $\| \cdot
\|$. In addition, we denote:
\[
\widetilde{\mu(A)} = \min(\mu(A),1-\mu(A))
\]
for all $A \subset X$. Lastly, we denote the Lebesgue measure on
$\Real^n$ by $\mes_n$.

\subsection{Isoperimetric inequalities}

The first topic of this note is an isoperimetric inequality for
$\mu$. In the setting of Example~\ref{e2} (and actually in a much
more abstract one), Bakry and Ledoux proved \cite{BakryLedoux} the
following isoperimetric inequality:

\begin{thm*}[Bakry -- Ledoux]
If the measure $\mu$ satisfies (\ref{ulc}) with $\| \cdot \| = |
\cdot |$ and $\delta(t) = t^2/8$, then for any $A \subset \RR^n$:
\begin{equation}\label{bl}
\mu^+_{|\cdot|}(A) \geq \phi\brac{\Phi^{-1}\Big( \widetilde{\mu(A)}
\Big)}~.
\end{equation}
Here as usual $\phi(t) = \frac{1}{\sqrt{2\pi}} \exp(-t^2/2)$ and
$\Phi(t) = \int_{- \infty}^t \phi(s) \, ds$.
\end{thm*}

This theorem is a generalisation of the isoperimetric inequality for
the Gaussian measure, proved by Sudakov, Tsirelson, and Borell
\cite{SudakovTsirelson,Borell-GaussianIsoperimetry}. In
\cite{BobkovLocalizedProofOfGaussianIso}, Bobkov gave a proof of the
Bakry--Ledoux inequality using the localisation technique; the
latter was introduced by Gromov and Milman \cite{Gromov-Milman} and
developed by Kannan, Lov\'asz and Simonovits
\cite{LSLocalizationLemma}, \cite{KLS} (see also
Gromov \cite[\S $3\frac{1}{2}.27$]{Gromov}). We extend
Bobkov's approach to the general case (\ref{ulc}) and prove:

\begin{thm} \label{isop1}
Suppose $\mu$ satisfies (\ref{ulc}) and (\ref{eq:delta-assumption}).
Then
\begin{equation} \label{isop_gamma_ineq}
\mu^+_{\| \cdot \|} (A)
    \geq C_\delta \widetilde{\mu(A)} \,
        \gamma \left(\log \frac{1}{\widetilde{\mu(A)}}\right)
    \; \text{ for all } \; A \subset \RR^n~,
\end{equation}
where:
\[\begin{split}
C_\delta &= \frac{e-1}{2e\max(2 \delta(\int_0^{+\infty} \exp(-2
\delta(t)) dt), 1)}, \quad \gamma(t) =
\frac{t}{\delta^{-1}(t/2)}~, \\
&\text{and} \quad \widetilde{\mu(A)} = \min(\mu(A), 1-\mu(A))~.
\end{split}\]
\end{thm}

\begin{cor}\label{isop_p}
Let $\delta(t) = \alpha t^p$ for $p \geq 2$ and $\alpha>0$ in the
setting of the previous theorem. Then:
\begin{equation}\label{isop_p_ineq}
\mu^+_{\| \cdot \|} (A)
    \geq c \alpha^{1/p} \widetilde{\mu(A)} \, \log^{1-1/p} \frac{1}{\widetilde{\mu(A)}}~,
\end{equation}
where $c > 0$ is a universal constant (independent of $p$).
\end{cor}

\begin{rem}
Note that $p$ can not be less than $2$; this follows from a second-order
Taylor expansion of $g$ in (\ref{ulc}).
\end{rem}

\begin{rem}
For $p=2$, Corollary~\ref{isop_p} recovers the Bakry--Ledoux Theorem
up to a universal constant: indeed,
\[ \phi(\Phi^{-1}(t))
    \leq C' t \sqrt{\log 1/t}~, \quad 0 \leq t \leq 1/2~.\]
\end{rem}

\begin{rem}\label{rem::bref}
In \cite{BobkovGaussianIsoLogSobEquivalent}, Bobkov proved that the following
inequality holds for any log-concave measure $\mu$ and any $r > 0$:
\begin{multline}\label{bisop}
\mu_{\| \cdot \|}^+(A) \geq \frac{1}{2r} \Big\{ \mu(A) \log \frac{1}{\mu(A)} \\
        + (1-\mu(A)) \log \frac{1}{1 - \mu(A)}
        + \log \mu \left\{\| x \| \leq r \right\} \Big\}~.
\end{multline}
In particular, (\ref{bisop}) implies a non-trivial isoperimetric
inequality for measures satisfying
(\ref{ulc}--\ref{eq:delta-assumption}). However, this inequality
would become weaker in higher dimension, whereas our results are
dimension-free.
\end{rem}

\subsection{Application: Uniformly convex bodies}

As before, let $V = (\RR^n, \| \cdot \|)$ be a normed space. The
volume measure $\lambda = \lambda_V$ on the unit ball of $V$ is
defined by:
\begin{equation}\label{unif}
\lambda = \frac{{\mes_n}|_{\{ \|x\| \leq 1\}}}{\mes_n(\{ \|x\| \leq
1\})} ;
\end{equation}
it arises naturally in geometric applications.

We would like to prove an isoperimetric inequality for $\lambda$,
with respect to the norm $\| \cdot \|$. It is easy to see that
$\lambda$ never satisfies the condition (\ref{ulc}) with $\delta >
0$. Therefore we follow the approach introduced by Bobkov and Ledoux
\cite{BobkovLedoux} and define an auxiliary measure $\mu$ that
satisfies (\ref{ulc}).

\subsubsection{$p$-uniformly convex bodies}

Choose $p \geq 2$, and let $\mu$ be the measure with density:
\begin{equation} \label{eq:mu-measure}
\frac{\exp(- \|x\|^p)}{\Gamma(1+n/p) \, \mes_n(\{\|x\| \leq 1\})}
\end{equation}
with respect to the Lebesgue measure.

\begin{prop*}[Bobkov -- Ledoux]
There exists a map $S: V \to V$ such that $S_\ast \mu = \lambda$ and
$\| S \|_\text{Lip} \leq C \left(\Gamma(1+n/p)\right)^{-1/n}$, where
$C>0$ is a universal constant.
\end{prop*}

It is clear that Lipschitz maps preserve isoperimetric inequalities,
so we may first establish one for $\mu$. The condition (\ref{ulc}) for $\mu$,
with $\delta(t) = \alpha t^p$, reads as
\begin{equation}\label{eq:puc_norm}
\frac{\|x\|^p+\|y\|^p}{2} - \left\|\frac{x+y}{2}\right\|^p
    \geq \alpha\|x - y\|^p \text{ for all } x,y \in \Real^n.
\end{equation}
This is one of the definitions of a $p$-uniformly convex norm (cf.\
Pisier \cite{PisierUniformlyConvexThm}).

\begin{ex}
The $\ell_q$ norm $\| \cdot \|_q$, $1<q<\infty$, satisfies
(\ref{eq:puc_norm}) with
\[ p = \begin{cases} 2, &q < 2 \\
                     q, &q \geq 2 \end{cases}~, \quad
\alpha = \begin{cases} \frac{q-1}{4}, &q < 2 \\
                     2^{-q}, &q \geq 2 \end{cases} .~\]
In fact, the same estimates holds for the space $L_q$. The case $q
\geq 2$ is due to Clarkson \cite{Clarkson} (see also Hanner
\cite{Hanner}), while the case $q<2$ follows from an unpublished
argument of Ball and Pisier (see Ball, Carlen and Lieb
\cite{BallCarlenLieb}).
\end{ex}

Therefore, if $\| \cdot \|$ is $p$-uniformly convex with coefficient
$\alpha$ (that is, if (\ref{eq:puc_norm}) holds), we can apply
Corollary~\ref{isop_p} and deduce (\ref{isop_p_ineq}). Combining
with the Bobkov -- Ledoux proposition above, we obtain the
following:

\begin{thm} \label{isop2}
Suppose the space $V$ is $p$-uniformly convex with constant $\alpha$
(that is, satisfies (\ref{eq:puc_norm})); let $\lambda$ be the
uniform measure on the unit ball of $\| \cdot \|$ (as in
(\ref{unif})). Then for any $A \subset \RR^n$:
\begin{equation}\label{eq:isop2}
\lambda^+_{\| \cdot \|} (A)
    \geq C \alpha^{1/p} n^{1/p} \widetilde{\lambda(A)} \,
        \log^{1-1/p} \frac{1}{\widetilde{\lambda(A)}}~,
\end{equation}
where $C > 0$ is a universal constant.
\end{thm}

This theorem continues the study of iso\-pe\-rimet\-ric properties
of $p$-uni\-for\-mly convex bodies by Bobkov and Zegarlinski
\cite[Ch.~14]{BobkovZegarlinski}. In particular, when $\lambda(A)$
is not exponentially small in the dimension, the inequality in
Theorem~\ref{isop2} improves the bound in \cite[Theorem 14.6]{BobkovZegarlinski}.
Under the same restriction, (\ref{eq:isop2}) improves (\ref{bisop})
with $r=1$ (which is however best possible in the class of all convex bodies).

\begin{rem} 
Here, as well as in Theorem~\ref{isop3} below, one may use an isoperimetric
inequality due to Barthe \cite{BartheGAFA} (which extends (\ref{bisop}))
and get a better bound for exponentially small sets. We do not pursue this point.
\end{rem}

\subsubsection{General uniformly convex bodies} \label{sec:intro-uniformly-convex-bodies}

We also generalise the above results to arbitrary uniformly convex
spaces. Recall that the modulus of convexity
$\delta_V : [0,2] \rightarrow [0,1]$
of a normed space $V = (X,\norm{\cdot})$ is defined as:
\[
\delta_V(\eps) = \inf \set{ 1- \norm{\frac{x+y}{2}} ;
\norm{x},\norm{y} \leq 1 \;,\; \norm{x-y} \geq \eps }.
\]
The space is called uniformly convex if $\delta_V(\eps)>0$ for all
$\eps>0$. From the works of Figiel \cite{FigielModulusOfConvexity},
Figiel-Pisier \cite{FigielPisier} and Pisier
\cite{PisierUniformlyConvexThm}, it is known that if:
\begin{equation} \label{eq:puc2_norm}
\delta_V(\eps) \geq \alpha' \eps^p \text{ for all } \eps \in [0,2],
\end{equation}
then (\ref{eq:puc_norm}) holds with $\alpha = \min(c,\alpha' /
2^p)$, and that if (\ref{eq:puc_norm}) holds then
(\ref{eq:puc2_norm}) holds with $\alpha' = \alpha / p$ (here $c>0$
is a universal constant). A space is therefore $p$-uniformly convex
if either (\ref{eq:puc_norm}) or (\ref{eq:puc2_norm}) hold, it is
however important to specify which definition one uses if the
dependence on $p$ is of interest.

In Section~\ref{sec:gucb} we derive the following proposition
from the results of Figiel-Pisier \cite{FigielPisier}:

\begin{prop} \label{prop:FigielPisier}
For all $x,y \in X$ such that $\norm{x}^2 + \norm{y}^2 \leq 2$, one
has:
\[
\frac{\norm{x}^2 + \norm{y}^2}{2} - \norm{\frac{x+y}{2}}^2 \geq c \,
\delta_V\brac{\frac{\norm{x-y}}{4}},
\]
where $c>0$ is a universal constant.
\end{prop}

Returning to the case $X = \Real^n$, choose $\mu$ to be the
probability measure with density:

\begin{equation} \label{eq:gen-density}
f(x) = \frac{1}{Z} \exp\brac{- \frac{n}{c} \norm{4 x}^2}
\mathbf{1}\left\{\norm{x}\leq \frac{1}{4}\right\}
\end{equation}
with respect to the Lebesgue measure, where $Z>0$ is a scaling
factor. Proposition \ref{prop:FigielPisier} clearly implies that
$\mu$ is uniformly log-concave, so we can apply Theorem \ref{isop1}
and deduce an isoperimetric inequality for $\mu$. To transfer this
inequality to the measure $\lambda_V$, we need to extend the
Bobkov--Ledoux proposition of the previous subsection. Our next
observation, which may be of independent interest, does precisely that.

\begin{dfn*}
A map $T: \RR^n \to \RR^n$ is called radial if it maps every ray to
itself in a monotone way; that is, if for every $x \neq 0$
\[\begin{cases}
T(\RR_+x) \subset \RR_+x &\text{and} \\
T|_{\RR_+ x}: \RR_+x \to \RR_+ x &\text{preserves the order on $\RR_+x$.}
\end{cases}\]
\end{dfn*}

Let $d\mu = f d\mes_n$ be an even log-concave probability measure (with
log-concave density $f$). Denote
\begin{equation}\label{def:kf}
K_f = \left\{ x \in \RR^n;
    \, n \int_0^{+\infty} f(rx) r^{n-1} dr \geq 1 \right\}~;
\end{equation}

It is not hard to see (cf.\ Proposition~\ref{prop:radial}) that there exists
a canonical radial map $T_f$ pushing forward $\mu$ to the restriction $\lambda$
of the Lebesgue measure to $K_f$.

K.~Ball showed \cite{Ball-kdim-sections} that $K_f$ is a symmetric convex
body; in other words, the unit ball of a norm $\norm{\cdot}_{K_f}$. In
Section~\ref{sec:Lip} we prove the following result (in a slightly more
general form):

\begin{thm} \label{thm:Lip}
Let $\mu = f d\mes_n$ be an even log-concave probability measure (with
log-concave density $f$); let $\lambda$ denote the restriction of the
Lebesgue measure on $K_f$, and let $T=T_f$ denote the canonical
radial map such that $T_* \mu = \lambda$. Then as a map $T: V \rightarrow V$
where $V = (\Real^n, \norm{\cdot}_{K_f})$, we have
$\norm{T}_{Lip} \leq C f(0)^{1/n}$, where $C>0$ is a universal
constant.
\end{thm}

\begin{rem}
The Bobkov-Ledoux proposition above is a particular case of the last
Theorem (up to another universal constant). We provide the details
at the end of Subsection~\ref{subsec:lip2}.
\end{rem}

In Section~\ref{sec:gucb} we apply Theorems~\ref{isop1} and \ref{thm:Lip}
to deduce the following:

\begin{thm} \label{isop3}
Let $V = (\Real^n,\norm{\cdot})$ be a uniformly convex space, and
let $\delta = \delta_V$ denote its modulus of convexity. Let
$\lambda = \lambda_V$ denote the uniform measure on the unit-ball of
$V$ (as in (\ref{unif})) and let $A \subset \RR^n$. Then:
\begin{equation*}\label{eq:isop3}
\lambda^+_{\| \cdot \|} (A)
    \geq c' C_{n,\delta} \frac{\widetilde{\lambda(A)} \log \frac{1}{\widetilde{\lambda(A)}}}
    {\delta^{-1}\brac{\frac{1}{2 n}\log \frac{1}{\widetilde{\lambda(A)}}}} ,
\end{equation*}
where:
\begin{equation} \label{eq:C-delta}
C_{n,\delta} = \frac{e-1}{2e\max(n \delta(\int_0^{1/4} \exp(-2 n
\delta(t)) dt), 1)},
\end{equation}
and $c'>0$ is a universal constant.
\end{thm}

Note that when $\delta(t) = \alpha t^p$ ($p \geq 2$), Theorem
\ref{isop3} recovers Theorem \ref{isop2} up to a universal constant.

\subsection{Connection to functional inequalities and concentration}\label{s_ineq}

In this subsection we study some corollaries of the isoperimetric inequalities
of the form (\ref{isop_gamma_ineq}) and (\ref{isop_p_ineq}).

\subsubsection{Concentration}

It is well-known that an isoperimetric inequality can be
equivalently rewritten in global form. It will be convenient to use
this in the following formulation (see Bobkov and Zegarlinski
\cite[p. 46]{BobkovZegarlinski} for an equivalent form):

\begin{prop}\label{propA}
Let $\mu$ be a probability measure on $\RR^n$ satisfying
\begin{equation} \label{isop_gamma_ineq1}
\mu^+_{\| \cdot \|} (A)
    \geq \widetilde{\mu(A)} \,
        \gamma \left(\log \frac{1}{\widetilde{\mu(A)}}\right)
\end{equation}
for every Borel set $A \subset \RR^n$ and some continuous function
$\gamma: [\log 2, +\infty) \to \RR_+$. Then for any Borel set $B
\subset \RR^n$ and any $\eps > 0$:
\begin{equation}\label{tmp1}
1 - \mu(B_{\eps,\| \cdot \|})
    \leq \exp \left( - h_{1-\mu(B)}^{-1}(\eps) \right)~,
\end{equation}
where:
\begin{equation}
h_a(x) = \int_{\log 1/a}^x \frac{dy}{\gamma(y)}~;
\end{equation}
for $y < \log 2$, $\gamma(y)$ should be interpreted as
$\gamma(\log \frac{1}{1 - \exp(-y)})$.

Conversely, if $\mu$ satisfies (\ref{tmp1}) for any Borel set
$B \subset \RR^n$, then (\ref{isop_gamma_ineq1}) holds.
\end{prop}

\begin{cor}\label{cor:1.9}
Let $\mu$ be a measure on $\RR^n$ such that for all $A \subset
\Real^n$:
\begin{equation}\label{isop_a}
\mu^+_{\| \cdot \|} (A)
    \geq \mathbf{c}_0 \widetilde{\mu(A)} \, \log^{1-1/p}
    \frac{1}{\widetilde{\mu(A)}}~.
\end{equation}
Then for every $B \subset \RR^n$, $\mu(B) \geq 1/2$, and every $\eps > 0$,
\begin{equation}
1 - \mu(B_{\eps, \| \cdot \|})
    \leq \exp \left\{ -
        \left[ \log^{1/p} \frac{1}{1 - \mu(B)} + \frac{\mathbf{c}_0\eps}{p}\right]^p
    \right\}~.
\end{equation}
\end{cor}

In Subsection~\ref{subsec:conc} we combine Proposition~\ref{propA}
and Corollary~\ref{cor:1.9} with the results of the previous subsections,
to deduce a concentration inequality for uniformly convex bodies. Then we
compare this inequality with the Gromov--Milman theorem \cite{Gromov-Milman}.

For completeness, we prove Proposition~\ref{propA} in Subsection~\ref{subsec:concpr}.

\subsubsection{Functional inequalities}

An isoperimetric inequality can be written in a functional
form; this was brought forth by Maz$'$ya, Federer, and Fleming
\cite{Mazya,FedererFleming} in the early 60's and later adapted by
Bobkov and Houdr\'e \cite{BobkovHoudre} to the context of
probability measures.

\begin{prop*}[Bobkov--Houdr\'e]
Let $\mu$ be a probability measure on a nor\-med space $(\RR^n, \| \cdot \|)$,
and let $I: [0, 1/2] \to \RR_+$ be an increasing
continuous function such that $I(0)=0$. The following are
equivalent:
\begin{enumerate}
\item For any Borel set $A \subset \RR^n$,
\begin{equation}
\mu^+_{\| \cdot \|} (A) \geq I(\widetilde{\mu(A)})~;
\end{equation}
\item For any locally Lipschitz function $F: \RR^n \to [0,1]$ such
that
\begin{equation}\label{techn}
\mu \{ F = 1 \} \geq t \in (0, \, 1/2) \quad \text{and} \quad
   \mu \{ F = 0 \} \geq 1/2~,
\end{equation}
we have:
\begin{equation}
\int \| \nabla F \|_\ast \, d\mu \geq I(t)~,
\end{equation}
where
\[ \|\nabla F \|_\ast = \limsup_{y \to x} \frac{|F(y)-F(x)|}{\|y-x\|}~.\]
\end{enumerate}
\end{prop*}

Let us focus on the case $I(t) = \mathbf{c}_0 t \log^{1/q} 1/t$, where
$1/q = 1 - 1/p$. We have the following:
\begin{prop}\label{propB}
Suppose a probability measure $\mu$ on $(\RR^n, \| \cdot \|)$ satisfies
\begin{equation}\label{isop_p_ineq1}
\mu^+_{\| \cdot \|} (A)
    \geq \mathbf{c}_0 \, \widetilde{\mu(A)} \,
        \log^{1/q}\frac{1}{\widetilde{\mu(A)}}
\end{equation}
for all $A \subset \RR^n$. Then:
\begin{enumerate}
\item\label{cap1}
For any locally Lipschitz function $F: \RR^n \to [0,1]$ satisfying
(\ref{techn}), we have:
\begin{equation}\label{1.26}
\int \| \nabla F \|_\ast \, d\mu \geq \mathbf{c}_0 t \log^{1/q}
1/t~;
\end{equation}
\item\label{cap2}
for any locally Lipschitz function $F: \RR^n \to [0, 1]$ satisfying
(\ref{techn}), we have:
\begin{equation}
\int \| \nabla F \|_\ast^q \, d\mu
    \geq c \, \mathbf{c}_0^q \, t \log 1/t~,
\end{equation}
where $c > 0$ is a universal constant;
\item \label{qls}
for any locally Lipschitz function $F: \RR^n \to \RR_+$,
\begin{equation}\label{qlogsob}
\int \| \nabla F \|_\ast^q d\mu
    \geq c' \mathbf{c}_0^q \int F^{q} \log \frac{F^q}{\int F^q d\mu} d\mu~,
\end{equation}
where $c' > 0$ is a universal constant.
\end{enumerate}
\end{prop}

Of course, \ref{cap1}.\ follows from the previous proposition (and
in fact, (\ref{1.26}) is equivalent to (\ref{isop_gamma_ineq1})).
Then, \ref{cap1}.\ implies \ref{cap2}.\ via standard arguments that
we reproduce for completeness in Section~\ref{sec:5}. Finally,
\ref{cap2}.\ is equivalent to \ref{qls}.\ (up to universal
constants); this is a reformulation of the arguments developed by
Bobkov and Zegarlinski \cite[Chapter 5.]{BobkovZegarlinski} in the
language of capacities put forth by Barthe and Roberto
\cite{BartheRoberto}.

The inequality (\ref{qlogsob}), called a $q$-log-Sobolev inequality,
was studied by Bobkov and Ledoux \cite{BobkovLedoux} and Bobkov and
Zegarlinski \cite{BobkovZegarlinski}. In particular, part
\ref{qls}.\ of the last proposition extends Theorem~16.3 in
\cite{BobkovZegarlinski}. Combining it with Theorems~\ref{isop1} and
\ref{isop2}, we recover the $q$-log-Sobolev inequalities proved by
Bobkov and Ledoux in \cite{BobkovLedoux}, up to universal constants.

\subsection{Acknowledgments.} We thank our supervisors Gideon Schechtman
and Vitali Milman for their guidance and support, and the referees for
careful reading. Part of this work was done while the authors
enjoyed the hospitality of the Henri Poincar\'e Institute in Paris.


\section{An Isoperimetric Inequality}\label{sec:isop}

\subsection{Reduction to one dimension}

This subsection is based on an argument that was introduced by
Gromov and Milman \cite{Gromov-Milman} to reduce the spherical
isoperimetric inequality to a certain one-dimensional fact; see also
Gromov \cite[\S $3\frac{1}{2}.27$]{Gromov} and Alesker
\cite{AleskerGM}. The corresponding argument in the affine case was
developed by Kannan, Lov\'asz and Simonovits
\cite{LSLocalizationLemma, KLS}, who also coined the term
`localisation lemma'; a different approach was put forth by
Fradelizi and Gu\'edon
\cite{FradeliziGuedonLocalization1,FradeliziGuedonLocalization2}.

We formulate the localisation lemma in terms of $\mu$-{\em needles}, as put
forth by S.~Bobkov; this corresponds to {\em convex descendants} in \cite{Gromov}.
It will be natural to work in an $n$-dimensional affine space $\mathbb{A}^n$
(cf.\ Remark~\ref{rem:aff}).

\medskip

Let $V = (\RR^n, \| \cdot \|)$ be a normed space acting by
translations on an affine space $\mathbb{A}^n$. Let $\mu$ be a
probability measure on $\mathbb{A}^n$ such that $\mu(H) = 0$ for
every affine hyperplane $H \subset \mathbb{A}^n$.

\begin{dfn*}
A (probability) measure $\sigma$ supported on an affine line $L
\subset \mathbb{A}^n$ (and not on any point) is called a
$\mu$-needle if
\[ \sigma = \lim_{k \to +\infty} \mu|_{C_k} / \mu(C_k) \]
is the weak limit of the scaled restrictions of $\mu$ to convex sets
\[ C_1 \supset C_2 \supset C_3 \cdots, \quad \mu(C_k)>0~.\]
\end{dfn*}

If the measure $\mu$ admits a lower semicontinuous density $f$ with
respect to the Lebesgue measure, the definition can be made more
explicit (see
\cite{LSLocalizationLemma,FradeliziGuedonLocalization2}). We will
only use the following property (see e.g.\
\cite[Lemma~2.5]{LSLocalizationLemma}):

\begin{prop1}{Description of $\mu$-needles}
If $\nu$ is a $\mu$-needle supported on $L$, then $\mu$ is absolutely
continuous with respect to the Lebesgue measure on $L$, and its density is
equal to $f|_L\phi$ for some log-concave function $\phi$ on $L$.
\end{prop1}

\begin{prop1}{Localisation principle: global form}
Let $\mu$ be a probability measure on $\mathbb{A}^n$ such that
$\mu(H)=0$ for every affine hyperplane $H \subset \RR^n$; let $a , b
\in (0, 1)$, $\eps > 0$. If every $\mu$-needle $\sigma$ supported on
an affine line $L_\sigma$ satisfies:
\begin{equation}\label{tmp2gf}
\sigma(A'_\eps) \geq b \text{ for every } A' \subset L_\sigma \text{
such that } \sigma(A') = a~,
\end{equation}
then also:
\begin{equation}\label{tmp3gf}
\mu(A_\eps) \geq b \text{ for every } A \subset \mathbb{A}^n \text{
such that } \mu(A) = a~.
\end{equation}
\end{prop1}

This is essentially the first step in \cite{LSLocalizationLemma}. It
will be more convenient to obtain an infinitesimal form of this
localisation principle. Given an isoperimetric inequality in the
general form:
\[
\mu^+(A) \geq I(\widetilde{\mu(A)}),
\]
where $I: [0,1/2] \rightarrow \Real_+$ is a continuous function, we
may of course write $I(a) = a \gamma(\log 1/a)$ for some continuous
function $\gamma: \Real_+ \rightarrow \Real_+$, obtaining the form
in (\ref{isop_gamma_ineq1}). By Proposition~\ref{propA}, a local
isoperimetric inequality of the form (\ref{isop_gamma_ineq1}) is
equivalent to the global inequality (\ref{tmp1}). Applying this
twice, we deduce the following:

\begin{prop1}{Localisation principle: local form}
Let $\mu$ be a probability measure on $\mathbb{A}^n$ such that
$\mu(H) = 0$ for every affine hyperplane $H \subset \mathbb{A}^n$,
and let $I: [0,1/2] \rightarrow \Real_+$ denote a continuous
function. If every $\mu$-needle $\sigma$ supported on an affine line
$L_\sigma$ satisfies:
\begin{equation}\label{tmp2}
\sigma_{\norm{\cdot}}^+ (A') \geq I\brac{\widetilde{\sigma(A')}} \text{ for every } A'
\subset L_\sigma \, ,
\end{equation}
then also:
\begin{equation}\label{tmp3}
\mu_{\norm{\cdot}}^+ (A) \geq I\brac{\widetilde{\mu(A)}} \text{ for every } A \subset
\mathbb{A}^n .
\end{equation}
\end{prop1}

To complete the reduction to one dimension, let us show that ``if $\mu$ is
uniformly log-concave, its needles are also uniformly log-concave''. The following
lemma extends \cite{BobkovLocalizedProofOfGaussianIso},
\cite[\S $3 \frac{1}{2}.27$, Ex.\ (e)]{Gromov}.

\begin{lem}\label{lem:aff}
Let $V = (\RR^n, \norm{ \cdot })$ be a normed space acting by
translations on an affine space $\mathbb{A}^n$, and let $\delta:
\Real_+ \rightarrow \Real_+$. If a measure $\mu$ on $\mathbb{A}^n$
satisfies the uniform log-concavity condition (\ref{ulc}) with
respect to $\delta$ and $\norm{\cdot}$, then every $\mu$-needle
$\sigma$ supported on an affine line $L \subset \mathbb{A}^n$
satisfies (\ref{ulc}) with respect to $\delta$ and the restriction
of $\norm{ \cdot }$ to the tangent space $L - L$.
\end{lem}

\begin{proof}[Sketch of proof]
Let $f$ denote the density of $\mu$ with respect to the Lebesgue
measure on $\mathbb{A}^n$. $\mu$ satisfies (\ref{ulc}), hence $f$ is
in particular log-concave. The super-level sets of $f$ are convex, hence
$f$ is equivalent to a lower semi-continuous density. Therefore the
description of needles formulated above is valid.

Now, $f|_L$ satisfies (\ref{ulc}) with respect to $\delta$ and
$\| \cdot \|_{L-L}$. Since $\phi$ satisfies (\ref{ulc}) with respect to
$\delta' \equiv 0$, it follows that $f|_L \phi$ satisfies
(\ref{ulc}) with respect to $\delta + \delta' = \delta$ and $\|
\cdot \|_{L-L}$.
\end{proof}

By the lemma, it is sufficient to prove Theorem~\ref{isop1} for $n=1$. In this
case, we only need the following property of one-dimensional uniformly log-concave
measures:

\begin{lem} \label{lem:g-majoration}
Let $V = (\Real^n,\norm{\cdot})$, and assume that $g : V \rightarrow
\Real \cup \set{+\infty}$ satisfies (\ref{ulc}). Assume in addition
that $a$ is a minimum point of $g$. Then:
\begin{equation}\label{eq:tail}
g(x) - g(a) \geq 2 \delta(\norm{x-a}),
\end{equation}
for all $x \in \Real^n$.
\end{lem}
\begin{proof}
If $g(x) = +\infty$, the claim is trivial. Otherwise, apply
(\ref{ulc}) with $y = a$. Then:
\[
\delta(\norm{x-a}) \leq \frac{g(x) - g(a)}{2} + g(a) -
g\brac{\frac{x+a}{2}} \leq \frac{g(x)-g(a)}{2},
\]
where we used the fact that $a$ is a minimum point of $g$ in the
last inequality.
\end{proof}

We will prove the isoperimetric inequality for one-dimensional
measures $\mu$ with density $f = \exp(-g)$, where $g$ satisfies
(\ref{eq:tail}). Any norm on $\RR^1$ is Euclidean, hence without
loss of generality $\| \cdot \| = | \cdot |$. Therefore
Theorem~\ref{isop1} is reduced to the following proposition (note
the factor 2 that we drop between (\ref{eq:tail}) and
(\ref{eq:g-condition}) to simplify the notation).

\begin{prop} \label{prop:1D-isoperimetry}
Let $\sigma$ denote a probability measure on $\Real$ with density
$f$. Assume that $f = \exp(-g)$, where $g : \Real \rightarrow \Real
\cup \set{+\infty}$ is a convex function with minimum at 0 and such
that:
\begin{equation} \label{eq:g-condition}
g(x) - g(0) \geq \delta(|x|)
\end{equation}
for all $x \in \Real$, and $\delta :\Real_+ \rightarrow \Real_+ \cup
\set{+\infty}$ satisfies (\ref{eq:delta-assumption}). Then:
\begin{equation} \label{eq:1D-isoperimetry}
\sigma^+(A) \geq C_\delta \widetilde{\sigma(A)}
\gamma\brac{\log\frac{1}{\widetilde{\sigma(A)}}}
\end{equation}
for any $A \subset \Real$, where
\[
C_\delta = \frac{e-1}{2e\max(\delta(\int_0^{+\infty}
\exp(-\delta(t)) dt), 1)}, \quad \gamma(t) =
\frac{t}{\delta^{-1}(t)}~.
\]
\end{prop}

\subsection{Proof of the one-dimensional inequality}

Before proceeding to the proof of Proposition
\ref{prop:1D-isoperimetry}, we collect several easy observations,
using the same notation as in the proposition.

\begin{lem} \label{lem:increasing-range}
The function $\gamma$ is non-decreasing. The function $x
\gamma\brac{ \log \frac{1}{x} }$ is strictly increasing on
$[0,1/e]$.
\end{lem}
\begin{proof}
The first part follows since $\delta(x)/x$ is non-decreasing by our
assumption (\ref{eq:delta-assumption}). For the second part, write:
\[
x \gamma\brac{ \log \frac{1}{x} } = \frac{ x \log \frac{1}{x} }{
\delta^{-1}(\log \frac{1}{x} ) },
\]
so the claim follows since $\delta$ (and hence $\delta^{-1}$) is
non-decreasing, whereas $x \log \frac{1}{x} $ is increasing on
$[0,1/e]$.
\end{proof}

Now denote:
\[
M_\delta = \int_0^\infty \exp(-\delta(x)) \, dx~.
\]

\begin{lem} \label{lem:max-bound}
\[
\exp(-g(0)) \geq \brac{2 M_\delta}^{-1}~.
\]
\end{lem}
\begin{proof}
Using $\int f(x) dx = 1$ and (\ref{eq:g-condition}):
\[
1 = \int_\Real \exp(-g(x)) dx \leq \exp(-g(0)) \int_\Real
\exp(-\delta(|x|) dx.
\]
\end{proof}

\begin{lem} \label{lem:trivial-bound}
If $g$ is a convex function on $\Real$ with minimum at 0, then for
all $x>0$:
\[
\int_x^\infty \exp(-g(y)) dy \leq \frac{x}{g(x)-g(0)}\exp(-g(x)).
\]
\end{lem}
\begin{proof}
By convexity, it follows that for all $y \geq x$:
\[
g(y) \geq \frac{g(x) - g(0)}{x} (y-x) + g(x).
\]
Using this to bound $\int_x^\infty \exp(-g(y)) dy$ from above, the
claim follows.
\end{proof}

Given a finite measure $\mu$ on $\Real$, we denote by $m(\mu)$
its median, i.e.\ (any) number $m$ for which $\mu((-\infty,m])
\geq \mu(\RR)/2$ and $\mu([m,\infty)) \geq \mu(\RR)/2$.

\begin{lem} \label{lem:median-max}
For any finite log-concave measure $d\mu = f dx$ on $\RR$,
\begin{equation}\label{eq:triv}
f(m(\mu)) \geq \frac{1}{2} \max_{x \in \Real} f(x).
\end{equation}
\end{lem}
\begin{proof}
Without loss of generality, assume $m = m(\mu) > 0$, $f(0) = \max f$
and $f(m) < f(0)$. Then $f$ is non-increasing on $\RR_+$. Replace
$\mu$ with $\mu|_{\RR_+}$; then the left-hand side of
(\ref{eq:triv}) may only decrease, whereas the right-hand side
retains its value.

Now replace $f$ by a log-affine function $f_1$ on $\Real_+$ such
that $f_1(0) = f(0)$ and $f_1(m) = f(m)$. In other words $f_1(x) =
\exp(-ax+b)|_{\RR_+}$, and our assumptions imply that $a>0$. Setting
$d \mu_1 = f_1 dx$, $\mu_1$ is a finite measure. Then $f_1 \leq f$
on $[0, m]$ and $f_1 \geq f$ on $[m, +\infty)$; hence $m(\mu_1) \geq
m(\mu)$ and $f(m(\mu)) = f_1(m(\mu)) \geq f_1(m(\mu_1))$.

Finally,
\begin{equation*}
f_1(m(\mu_1)) = \frac{1}{2} \max_{x \in \Real_+} f_1(x);
\end{equation*}
this concludes the proof.
\end{proof}

\begin{proof}[Proof of Proposition \ref{prop:1D-isoperimetry}]
By a general result of Bobkov (\cite[Proposition
2.1]{BobkovExtremalHalfSpaces}) on extremal isoperimetric sets of
log-concave densities, it is enough to verify
(\ref{eq:1D-isoperimetry}) on sets $A$ of the form $(-\infty,a]$ and
$[b,\infty)$. Given a point $x \in \Real$, denote $A = [x,\infty)$
if $x \geq 0$ and $A = (-\infty,x]$ if $x < 0$. We will show that
the set $A$ satisfies:
\[
\sigma^+(A) \geq C_\delta \widetilde{\sigma(A)}
\gamma\brac{\log\frac{1}{\widetilde{\sigma(A)}}},
\]
and this will conclude the proof. Assume w.l.o.g. that $x \geq 0$,
since our hypotheses are symmetric about the origin.

First, recall that by another result of Bobkov
(\cite[Proposition 4.1]{BobkovGaussianIsoLogSobEquivalent}), a
log-concave probability measure $\mu$ with density $f$ on $\Real$
always satisfies the following Cheeger-type isoperimetric
inequality:
\[
\mu^+(A) \geq 2 f(m) \min(\mu(A),1-\mu(A)),
\]
where $m$ is the median of $\mu$. Together with Lemma
\ref{lem:median-max}, this implies:
\begin{equation} \label{eq:Cheeger}
\sigma^+(A) \geq \exp(-g(0))
\widetilde{\sigma(A)}.
\end{equation}
Loosely speaking, this Cheeger-type inequality will take care of the
case when $\widetilde{\sigma(A)}$ is large. The case when
$\widetilde{\sigma(A)}$ is small will be handled by Lemma
\ref{lem:trivial-bound}, which, together with the assumption
(\ref{eq:g-condition}) and the fact that $\delta$ is increasing,
imply that for any $x>0$:
\[
\sigma(A) = \int_x^\infty \exp(-g(y)) dy \leq
\frac{\delta^{-1}(g(x)-g(0))}{g(x)-g(0)} \exp(-g(x)).
\]
Recalling the definition of $\gamma$ and denoting
$\sigma^+_{\max} = \exp(-g(0))$, this means:
\begin{equation} \label{eq:whole-point}
\sigma(A) \leq \frac{\sigma^+(A)}{\gamma(g(x)-g(0))} =
\frac{\sigma^+(A)}{\gamma(\log \frac{\sigma^+_{\max}}{\sigma^+(A)})}.
\end{equation}
This inequality is almost what we need, and the
rest of the proof will be dedicated to replacing $\sigma^+$
with $\sigma$ inside the $\gamma$ function.

More formally, we distinguish between five cases.
\begin{enumerate}
\item
$\widetilde{\sigma(A)} \geq c_\delta$, where $c_\delta \leq 1/e$
depends solely on $\delta$ and will be determined later. In this
case, by (\ref{eq:Cheeger}) and Lemma \ref{lem:max-bound}:
\[
\sigma^+(A) \geq \exp(-g(0))
\widetilde{\sigma(A)} \geq \frac{1}{2 M_\delta }
\widetilde{\sigma(A)}.
\]
The function $\gamma$ is non-decreasing by Lemma \ref{lem:increasing-range},
therefore
\[
\sigma^+(A) \geq \frac{1}{2 M_\delta \; \gamma(\log \frac{1}{c_\delta})}
     \widetilde{\sigma(A)} \gamma(\log \frac{1}{\widetilde{\sigma(A)}}).
\]
\item \label{b}
$1 - \sigma(A) = \widetilde{\sigma(A)} < c_\delta$ and $g(x)-g(0) <
\log \frac{1}{c_\delta}$. Using (\ref{eq:g-condition}):
\[
\sigma(A) \leq  \int_0^\infty \exp(-g(y)) dy \leq \exp(-g(0))
\int_0^\infty \exp(-\delta(y)) dy,
\]
and since $g(x)-g(0) < \log \frac{1}{c_\delta}$ we conclude that:
\[
1 - c_\delta < \sigma(A) \leq \frac{1}{c_\delta} \exp(-g(x))
M_\delta = \frac{M_\delta}{c_\delta} \sigma^+(A).
\]
By Lemma \ref{lem:increasing-range}, $x \gamma( \log \frac{1}{x} )$
is monotone increasing on $[0,1/e]$. Since $\widetilde{\sigma(A)} <
c_\delta \leq 1/e$, we conclude that:
\[
\sigma^+(A) \geq \frac{(1-c_\delta) c_\delta \gamma(\log
\frac{1}{c_\delta}) }{ M_\delta \gamma(\log\frac{1}{c_\delta}) }
\geq \frac{(1-c_\delta)}{ M_\delta \gamma(\log\frac{1}{c_\delta}) }
\widetilde{\sigma(A)} \gamma(\log\frac{1}{\widetilde{\sigma(A)}}).
\]
\item
$\sigma(A) = \widetilde{\sigma(A)} < c_\delta$ and $g(x)-g(0) < \log
\frac{1}{c_\delta}$. As in \ref{b}.\ :
\begin{multline*}
1 - \sigma(A) = \int_{-\infty}^0 \exp(-g(y)) dy + \int_0^x
\exp(-g(y)) dy \\
\leq \exp(-g(0)) M_\delta + \exp(-g(0)) x \leq
\frac{1}{c_\delta} \exp(-g(x)) (M_\delta + x).
\end{multline*}
Using (\ref{eq:g-condition}) and the inequality
$g(x)-g(0) < \log \frac{1}{c_\delta}$,
\[
x \leq \delta^{-1}(g(x)-g(0)) \leq \delta^{-1}(\log
\frac{1}{c_\delta}).
\]
Hence:
\[
1-c_\delta \leq 1 - \sigma(A) \leq \frac{M_\delta + \delta^{-1}(\log
\frac{1}{c_\delta})}{c_\delta} \sigma^+(A).
\]
Now choose
\begin{equation} \label{eq:cd-definition}
c_\delta := \min(1/e,\exp(-\delta(M_\delta))),
\end{equation}
which yields:
\[
\sigma^+(A) \geq \frac{(1-c_\delta) c_\delta}{2 \delta^{-1}(\log
\frac{1}{c_\delta})} = \frac{(1-c_\delta) c_\delta \gamma(\log
\frac{1}{c_\delta})}{2\log \frac{1}{c_\delta}}.
\]
By the monotonicity of $x \gamma( \log \frac{1}{x} )$ as in
\ref{b}.\ , we conclude that:
\[
\sigma^+(A) \geq \frac{(1-c_\delta)}{2\log \frac{1}{c_\delta}}
\widetilde{\sigma(A)} \gamma(\log \frac{1}{\widetilde{\sigma(A)}}).
\]
\item
$\widetilde{\sigma(A)} < c_\delta$, $g(x)-g(0) \geq \log
\frac{1}{c_\delta}$ and $\frac{\sigma^+(A)}{\gamma(g(x)-g(0))} \geq
1/e$. Since $\gamma$ is non-decreasing:
\[
\sigma^+(A) \geq \frac{1}{e} \gamma(g(x)-g(0)) \geq \frac{1}{e
c_\delta} c_\delta \gamma( \log \frac{1}{c_\delta} ).
\]
Using the monotonicity of $x \gamma( \log \frac{1}{x} )$ as in
\ref{b}.\ , we conclude that:
\[
\sigma^+(A) \geq \frac{1}{e c_\delta} \widetilde{\sigma(A)} \gamma(
\log \frac{1}{\widetilde{\sigma(A)}} ).
\]
\item
$\widetilde{\sigma(A)} < c_\delta$, $g(x)-g(0) \geq \log
\frac{1}{c_\delta}$ and $\frac{\sigma^+(A)}{\gamma(g(x)-g(0))} <
1/e$. Recall that by (\ref{eq:whole-point}):
\[
\sigma(A) \leq \frac{\sigma^+(A)}{\gamma(g(x)-g(0))} < \frac{1}{e},
\]
implying in particular that $\widetilde{\sigma(A)} = \sigma(A)$. We
will show:
\begin{equation} \label{e}
\sigma^+(A) \geq D_\delta \frac{\sigma^+(A)}{\gamma(g(x)-g(0))}
\gamma\brac{ \log \frac{\gamma(g(x)-g(0))}{\sigma^+(A)} },
\end{equation}
which by the monotonicity of $x \gamma( \log \frac{1}{x} )$ on
[0,1/e] will imply:
\begin{equation} \label{eq:e-needed}
\sigma^+(A) \geq D_\delta \widetilde{\sigma(A)} \gamma( \log
\frac{1}{\widetilde{\sigma(A)}} ).
\end{equation}
Denote $V_x = g(x) - g(0)$. Then (\ref{e}) is equivalent to showing:
\[
\frac{\gamma(V_x (1 + \frac{\log
\frac{\gamma(V_x)}{\exp(-g(0))}}{V_x}))}{\gamma(V_x)} \leq
1/D_\delta.
\]
Recall that $\gamma$ is non-decreasing and note that
$\frac{\gamma(x)}{x} = \frac{1}{\delta^{-1}(x)}$ is non-increasing.
Requiring that $D_\delta \leq 1$, it is therefore enough to show:
\[
1 + \frac{\log \frac{\gamma(V_x)}{\exp(-g(0))}}{V_x} \leq 1 /
D_\delta.
\]
Denoting $B_\delta := 1/D_\delta - 1$, the latter is equivalent to:
\[
\gamma(V_x) \leq \exp(B_\delta V_x) \exp(-g(0)),
\]
which from the definition of $\gamma$ is equivalent to:
\[
\delta(V_x \exp(-B_\delta V_x) \exp(g(0))) \leq V_x.
\]
The maximum of the function $z \mapsto z \exp(-B_\delta z)$ is equal to $1/(e
B_\delta)$, hence it is enough to require that:
\[
\delta\brac{\frac{\exp(g(0))}{e B_\delta}} \leq V_x.
\]
We have assumed that $V_x = g(x) - g(0) \geq \log\frac{1}{c_\delta}$; therefore
by the definition (\ref{eq:cd-definition}) of $c_\delta$ the following condition
will suffice:
\begin{equation}\label{star}
\frac{\exp(g(0))}{e B_\delta} \leq M_\delta.
\end{equation}
By Lemma \ref{lem:max-bound}, (\ref{star}) holds for $B_\delta = 2/e$
(independent of $\delta$ in fact!). To conclude, (\ref{eq:e-needed})
is satisfied with $D_\delta = \frac{e}{e+2}$.
\end{enumerate}

Summing up all the five requirements for the constant $C_\delta$ in the
conclusion of the proposition, we see that we can choose:
\[
C_\delta \leq \min \brac{\frac{1}{2 M_\delta \; \gamma(\log
\frac{1}{c_\delta})},\frac{(1-c_\delta)}{ M_\delta \gamma(\log
\frac{1}{c_\delta}) },\frac{(1-c_\delta)}{2\log \frac{1}{c_\delta}},
\frac{1}{e c_\delta},\frac{e}{e+2}}.
\]
From the definition (\ref{eq:cd-definition}) of $c_\delta$, we see
that $\log \frac{1}{c_\delta} = \max(\delta(M_\delta),1)$ and that
$\gamma(\log \frac{1}{c_\delta}) \leq \max(\delta(M_\delta),1) /
M_\delta$. It is then not hard to check that we can choose:
\[
C_\delta := \frac{e-1}{2 e \max(\delta(M_\delta),1)},
\]
as claimed.
\end{proof}

\subsection{A simpler proof with further assumptions}

Note that the uniform convexity (\ref{ulc}) of $g$ was not used in
the statement and proof of Proposition \ref{prop:1D-isoperimetry}.
We remark here that by using this property, we obtain a simpler
proof of a one-dimensional isoperimetric inequality, which may be
used to complete the proof of Theorem \ref{isop1} in place of
Proposition \ref{prop:1D-isoperimetry}. The key observation is the
following:

\begin{lem}\label{lemma:n1}
Suppose $g: (\RR,\abs{\cdot}) \to \RR \cup \{ +\infty\}$ satisfies
(\ref{ulc}), i.e.:
\begin{equation}\label{eq:ulc.1}
\frac{g(x)+g(y)}{2} - g \left( \frac{x+y}{2} \right) \geq
\delta(|x-y|) \geq 0~,
    \quad x,y \in \RR~.
\end{equation}
Then for any $x_0 \in \RR$ :
\begin{equation}\label{eq:decay}
g(x) \geq g(x_0) + g'(x_0) (x-x_0) + 2\delta(|x-x_0|)~,
\end{equation}
where $g'(x_0)$ is any value between $g'_l(x_0)$ and $g'_r(x_0)$,
the left and right derivatives at $x_0$, respectively.
\end{lem}

\begin{proof}
Immediate by applying Lemma \ref{lem:g-majoration} to the function
$g - g'(x_0) (x-x_0)$, which attains its minimum at $x_0$.
\end{proof}

\begin{prop}\label{prop:new}
Let $\sigma$ be a probability measure on $\RR$ such that
\[ d\sigma(x) = \exp(-g(x)) dx~, \]
where $g$ satisfies (\ref{eq:ulc.1}). Then
\begin{equation}\label{eq:isop.new}
\sigma^+(A)
    \geq \widetilde{\sigma(A)}
        \psi^{-1} \left( \frac{1}{2 \widetilde{\sigma(A)}}\right)~,
    \quad A \subset \RR~,
\end{equation}
where
\begin{equation}\label{eq:def.psi}
\psi(t) = t \phi(t)~, \quad \phi(t) = \int_0^{+\infty} \exp(t x -
2\delta(x)) dx~.
\end{equation}
\end{prop}

\begin{proof}
As before, by a general result of Bobkov (\cite[Proposition
2.1]{BobkovExtremalHalfSpaces}) on extremal isoperimetric sets of
log-concave densities, it is enough to verify (\ref{eq:isop.new}) on
sets $A$ of the form $(-\infty,x_0]$ and $[x_0,\infty)$. By
symmetry, we may restrict ourselves to sets $[x_0, +\infty)$,
$\sigma([x_0, \infty)) = a \leq 1/2$.

Denote $a^+ = \exp(-g(x_0))$. By (\ref{eq:decay}),
\begin{equation} \label{eq:new1}
a = \int_{x_0}^\infty \exp(-g(x)) dx \leq a^+ \phi(-g'(x_0))~,
\end{equation}
and similarly
\begin{equation} \label{eq:new2}
1/2 \leq 1 - a \leq a^+ \phi(g'(x_0))~.
\end{equation}

Now consider two cases.
\begin{description}
\item[Case 1: $g'(x_0) > 0$.]
By (\ref{eq:def.psi}), $\phi(-g'(x_0)) \leq 1/g'(x_0)$; hence
$g'(x_0) \leq a^+/a$ using (\ref{eq:new1}) and $a^+ \phi(a^+/a) \geq
a^+ \phi(g'(x_0)) \geq 1/2$ using (\ref{eq:new2}). Therefore
\[\psi(a^+/a) = (a^+/a) \phi(a^+/a) \geq 1/2a~,\]
which implies (\ref{eq:isop.new}).
\item[Case 2: $g'(x_0) \leq 0$.]
By (\ref{eq:new2}), $a^+ \phi(0) \geq a^+ \phi(g'(x_0)) \geq 1/2$,
hence
\begin{equation} \label{new3}
a^+ \geq \frac{1}{2\phi(0)}~.
\end{equation}
Next, since $\phi$ is monotone, $\phi(\frac{1}{2a \phi(0)}) \geq
\phi(0)$, hence $\psi \left( \frac{1}{2a \phi(0)} \right) \geq
\frac{1}{2a}$, and we conclude by (\ref{new3}) that:
\[
a^+ \geq  \frac{1}{2\phi(0)} \geq a \psi^{-1} \left( \frac{1}{2a}
\right)~.
\]
\end{description}
\end{proof}

\begin{rem}\label{rem:phi}
It is easy to verify that the function $\phi$ defined in
(\ref{eq:def.psi}) is log-convex, i.e. $\log \phi$ is convex.
\end{rem}

\begin{rem}
Note that when $\delta(t) = c t^p$ ($p \geq 2$), the inequalities
obtained in Propositions \ref{prop:1D-isoperimetry} and
\ref{prop:new} are equivalent, up to universal constants.
\end{rem}


\section{Lipschitz Maps} \label{sec:Lip}

This section is dedicated to the proof of an extended form of
Theorem \ref{thm:Lip}.

\begin{prop}\label{prop:radial}
Let $\mu$ be a finite absolutely continuous measure on $\RR^n$.
There exists a $\mu$-a.e.\ unique radial map $T$ that pushes $\mu$
forward to the restriction of the Lebesgue measure to some
star-shaped set $K \subset \RR^n$.

If $d\mu = f \, d\mes_n$, we may choose $K=K_f$ and $T=T_f$, where:
\begin{eqnarray}
\label{eq:def-K_f} K_f &=& \left\{ x \in \RR^n;
    \, v(x) \leq 1 \right\}~,\\
\nonumber v(x) &=& \brac{ n \int_0^{+\infty} f(rx) r^{n-1} dx
}^{-\frac{1}{n}}~,
\end{eqnarray}
and $T_f$ is given by $T_f(0)=0$ and:
\begin{equation} \label{eq:T_f}
 T_f(x) = \brac{\frac{\int_0^{1} f(r x) r^{n-1} dr}
                {\int_0^{\infty} f(r x) r^{n-1} dr}}^{\frac{1}{n}}
        \frac{x}{v(x)} \;\; , \;\; x\neq 0. \\
\end{equation}
\end{prop}

\begin{proof}[Proof of Proposition~\ref{prop:radial}]
Let $T: \RR^n \to \RR^n$ be a radial map pushing $\mu$ forward to the Lebesgue
measure restricted to a star-shaped body $K$. Define:
\[ w(x) = \inf \left\{ t > 0 \, ; \,  t^{-1} x \in K \right\}~;\]
then the restriction of $T$ to a ray $\RR_+ x$, $w(x) = 1$, has the form:
\[
rx \mapsto u(x, r) x~, \quad r > 0~.
\]
Passing to polar coordinates and using the Fubini theorem, we see
that $T_\ast \mu$ is equal to the restriction of $\mes_n$ to $K$
iff, for almost every ray $\RR_+ x$, $w(x) = 1$, the map
$u(x,\cdot)$ pushes $f(rx) r^{n-1} dr$ forward to $\mathbf{1}_{[0,
1]} r^{n-1} dr$; that is, if
\begin{equation}\label{my3.2}
\int_{0}^{1} \phi(r) r^{n-1} dr
    = \int_{0}^{\infty} \phi(u(x, r)) f(rx) r^{n-1} dr
\end{equation}
for any test function $\phi \in C_0(\RR_+)$.
Setting $\phi = \mathbf{1}_{[0, T]}$ in (\ref{my3.2}) and letting $T \to \infty$,
we see that \begin{equation}\label{my3.2'}
\frac{1}{n} = \int_0^{\infty} f(rx) r^{n-1} dr~.
\end{equation}
Hence $v(x) = 1$ for (almost) every $x$ such that $w(x) = 1$. Both $v$ and $w$
are homogeneous functions, hence $v(x) = w(x)$ for $\mu$-a.e.\ $x \in \RR^n$.

Now use $\phi=\mathbf{1}_{[0,u(x,s)]}$ in (\ref{my3.2}). Since
$u(x,\cdot)$ is monotone, we deduce:
\[ u(x, s)^n = n \int_{0}^s f(rx) r^{n-1} dr
    = \frac{\int_{0}^s f(rx) r^{n-1} dr}{\int_{0}^{\infty} f(rx) r^{n-1} dr}~, \]
at every point of continuity $s$ of $u(x,\cdot)$. Therefore
\begin{equation}\label{my3.2'''}
T(sx)
    = \brac{\frac{\int_{0}^s f(rx) r^{n-1} dr}
           {\int_{0}^{\infty} f(rx) r^{n-1} dr}}^{1/n}
      x, \quad v(x) = 1~,
\end{equation}
which is equivalent to (\ref{eq:T_f}).
\end{proof}
\begin{rem}
Note that in particular, $\mes_n(K_f) = \mes_n(K) = \mu(\Real^n)$.
\end{rem}

The following proposition was proved by K.~Ball
\cite{Ball-kdim-sections} for even log-concave functions and
extended by Klartag \cite[Theorem 2.2]{KlartagPerturbationsWithBoundedLK}
to general log-concave functions.

\begin{prop*}[Ball]
If $f$ is a log-concave function on $\RR^n$, then $K_f$ is a convex
body.
\end{prop*}

Note that we do not assume at this stage that $f$ is even. Therefore
$K_f$ may not necessarily be symmetric about the origin, so formally
we can not identify it with the unit-ball of some norm
$\norm{\cdot}_{K_f}$. Nevertheless, we denote:
\begin{equation} \label{eq:K_f-def}
\norm{x}_{K_f} = \left( n \int_0^{\infty} f(rx) r^{n-1} dr
\right)^{-\frac{1}{n}}~;
\end{equation}
by the above proposition, this is a convex function on $\RR^n$,
which is in addition homogeneous. By definition (\ref{eq:def-K_f}),
we have:
\[
K_f = \set{x \in \Real^n ; \norm{x}_{K_f} \leq 1}.
\]
In addition, we denote:
\[
\widehat{K_f} = K_f \cap -K_f
\]
which is now a convex body symmetric about the origin, and we
associate with it the corresponding norm
$\norm{\cdot}_{\widehat{K_f}}$.

\medskip
We can now state the following result, which extends Theorem
\ref{thm:Lip}:

\begin{thm}\label{thm:Lip1}
Let $f$ denote a log-concave function on $\Real^n$ with barycenter
at the origin such that $0 < \int f(x) dx < \infty$. Let $\mu$
denote the measure with density $f$, and let $\lambda$ denote the
restriction of the Lebesgue measure to $K_f$. Denote by $T=T_f$ the
canonical radial map (given by (\ref{eq:T_f})) such that $T_* \mu =
\lambda$, and let $u : (\Real^n,\norm{\cdot}_{\widehat{K_f}})
\rightarrow [0,1]$ be defined by:
\[
T(x) = u(x) \frac{x}{\norm{x}_{K_f}}
\]
for $x \neq 0$ and $u(0)=0$. Then $\norm{u}_{Lip} \leq C
f(0)^{1/n}$, where $C>0$ is a universal constant.
\end{thm}

When $f$ is in addition even, $\widehat{K_f} = K_f$ and
$\norm{\cdot}_{K_f}$ is indeed a norm. Theorem \ref{thm:Lip} is then
deduced from Theorem \ref{thm:Lip1} using the following lemma, which
was essentially proved by Bobkov and Ledoux \cite{BobkovLedoux}.

\begin{lem} \label{lem:T2u}
Let $V = (X,\norm{\cdot})$ denote a normed space, and let $T : V
\rightarrow V$ be the map defined by $T(0)=0$ and:
\[
T(x) = u(x) \frac{x}{\norm{x}}
\]
for $x\neq 0$, where $u: X \rightarrow \Real_+$ has a finite
Lipschitz constant and satisfies $u(0) = 0$. Then:
\[
\norm{T}_{Lip} \leq 3 \norm{u}_{Lip}.
\]
\end{lem}

\begin{proof}
Let $x,y \in X$. By continuity, we may assume that $x,y \neq 0$.
Then:
\[\begin{split}
&\norm{T(x) - T(y)} = \norm{u(x) \frac{x}{\norm{x}} -
u(y) \frac{y}{\norm{y}}} \\
&\qquad\leq \norm{u(x) \frac{x}{\norm{x}} - u(x) \frac{y}{\norm{y}}}
+ \norm{u(x) \frac{y}{\norm{y}} - u(y)
\frac{y}{\norm{y}}} \\
&\qquad= \abs{u(x) - u(0)} \norm{\frac{x}{\norm{x}} -
\frac{y}{\norm{y}}} + \abs{u(x)-u(y)} \\
&\qquad \leq \norm{u}_{Lip} \norm{x} \brac{\norm{\frac{x}{\norm{x}}
- \frac{y}{\norm{x}}} +
\norm{\frac{y}{\norm{x}} - \frac{y}{\norm{y}} }} + \norm{u}_{Lip}\norm{x-y} \\
&\qquad = \norm{u}_{Lip} \norm{x} \norm{y} \abs{\frac{1}{\norm{x}} -
\frac{1}{\norm{y}}} + 2 \norm{u}_{Lip}\norm{x-y} \leq 3
\norm{u}_{Lip} \norm{x-y}.
\end{split}\]
\end{proof}

For the proof of Theorem \ref{thm:Lip1}, we need to compile
several known results about log-concave functions.

\subsection{Additional Preliminaries}

Another convex body associated to a log-concave function $f$ on
$\Real^n$ was put forth by B. Klartag and V. Milman
\cite{KlartagMilmanLogConcave}. Assume that $f(0)>0$, we define the
(convex) body $K_f^0$ as the set:
\begin{equation} \label{eq:def-K_f^0}
K_f^0 = \set{x \in \Real^n ; f(x) \geq f(0) \exp(-n)}.
\end{equation}
We will use a relation between $K_f$ and $K_f^0$ that was proved
(under slightly different assumptions) by Klartag and Milman
\cite[Lemmata 2.1,2.2]{KlartagMilmanLogConcave}:

\begin{prop}[Klartag--Milman]\label{prop:K_f^0-relation}
Let $f$ be a log-concave density on
$\RR^n$, and assume that $f(0)>0$. Then:
\[ K_f \subset C_n (\sup_x f(x))^{\frac{1}{n}} K_f^0~,\]
where $C_n>1$ and $C_n \rightarrow 1$ as $n \rightarrow \infty$.
Moreover, if $f$ attains its maximum at $0$, then:
\[ f(0)^{\frac{1}{n}} K_f^0 \subset D_n K_f~, \]
where $D_n > 2$ and $D_n \rightarrow 2$ as $n \rightarrow \infty$.
\end{prop}

The next lemma is a one dimensional computation for log-concave
functions. For even functions, this fact goes back to
Ball \cite{Ball-PhD}, and Milman and Pajor \cite{Milman-Pajor-LK}.
For arbitrary log-concave functions, this was extended by Klartag
\cite[Lemma 2.6]{KlartagPerturbationsWithBoundedLK} as follows:

\begin{lem} \label{lem:moments-lemma}
Let $f: \Real_+ \rightarrow \Real_+$ denote a non-constant
log-concave function, and let $n \geq 1$. Assume that $f(0)=1$ and
that:
\begin{equation}\label{eq:sup}
\sup_x f(x) \leq \exp(n)~.
\end{equation}
Then:
\[
C_1 \leq \frac{n^{\frac{n+1}{n}}}{e (n+1)} \leq \frac{\int_0^\infty
f(r) r^{n} dr}{(\int_0^\infty f(r) r^{n-1} dr)^{\frac{n+1}{n}}} \leq
\frac{n!}{((n-1)!)^{\frac{n+1}{n}}} \leq C_2,
\]
where $C_1,C_2>0$ are universal constants. In fact, the assumption
(\ref{eq:sup}) is not needed for the right-hand side of the
inequality.
\end{lem}

The last proposition we need is due to M.~Fradelizi \cite[Theorem 4]{FradeliziCentroid}:
\begin{prop}[Fradelizi] \label{prop:Fradelizi}
Let $f$ denote a log-concave density on $\RR^n$
such that $0< \int f(x) dx < + \infty$, and let $x_0$ denote its barycenter.
Then:
\[
g(x_0) \geq \exp(-n) \sup_{x \in \Real^n} g(x).
\]
\end{prop}

\subsection{Proof of Theorem \ref{thm:Lip1}}\label{subsec:lip2}

By (\ref{eq:T_f}), $T(x) = u(x) \frac{x}{\norm{x}_{K_f}}$ for $x
\neq 0$, where $u$ is given by:
\begin{equation} \label{eq:u-definition}
u(x) = \brac{\frac{\int_0^1 r^{n-1} f(r x) dr}{ \int_0^{\infty}
r^{n-1} f(r x) dr}}^{\frac{1}{n}}
\end{equation}
for $x \neq 0$ and $u(0)=0$. We thus verify that $u$ is continuous
at 0.

\medskip
\noindent{\bf Step 1}: Reduction to smooth $f$.
\medskip

Define, for $\eps>0$, $f_\eps := f \ast \eps^{-n} G(x/\eps)$, where
$G$ is the standard Gaussian density on $\Real^n$ and $\ast$ denotes
convolution. Clearly $f_\eps$ is a smooth function with barycenter
at 0. By the Pr\'ekopa-Leindler Theorem, $f_\eps$ is log-concave, as
the convolution of two log-concave functions.

Let $\mu_\eps$ denote the measure with density $f_\eps$,
$\lambda_\eps$ the Lebesgue measure on $K_{f_\eps}$, and let
$T_\eps$ denote the map radially pushing forward the measure
$\mu_\eps$ onto $\lambda_\eps$. Let $u_\eps$ be defined by
\[
T_\eps(x) = u_\eps(x) \frac{x}{\norm{x}_{K_{f_\eps}}},
\]
with $u_\eps(0) = 0$. Given $x,y \in \Real^n$, it is clear from
(\ref{eq:u-definition}) and (\ref{eq:K_f-def}) that $u_\eps(x)
\rightarrow u(x)$, $u_\eps(y) \rightarrow u(y)$,
$\norm{x-y}_{K_{f_\eps}} \rightarrow \norm{x-y}_{K_f}$ and
$\norm{x-y}_{\widehat{K_{f_\eps}}} \rightarrow
\norm{x-y}_{\widehat{K_f}}$ as $\eps$ tends to 0. If we assume that
$\norm{u_\eps}_{Lip} \leq C f_\eps(0)^{1/n}$, we have:
\[
\abs{u_\eps(x) - u_\eps(y)} \leq C f_\eps(0)^{1/n} \norm{x -
y}_{\widehat{K_{f_\eps}}}.
\]
Passing to the limit as $\eps \rightarrow 0$, it follows that:
\[
\abs{u(x) - u(y)} \leq C f(0)^{1/n} \norm{x - y}_{\widehat{K_{f}}},
\]
and we conclude that $\norm{u}_{Lip} \leq C f(0)^{1/n}$. It is
therefore enough to restrict our discussion to smooth functions.

\medskip
\noindent{\bf Step 2}: Proof for smooth functions with $f(0)=1$.
\medskip

Assume that $f(0)=1$.

Note that since $f$ and thus $u$ are assumed to be smooth,
\[
\norm{u}_{Lip} = \sup_{x \in \Real^n} \norm{\nabla
u(x)}_{\widehat{K_f}}^\ast~,
\]
where $\norm{\cdot}_{\widehat{K_f}}^\ast = \sup_{h \in
\widehat{K_f}} \scalar{\cdot,h}$ is the dual norm to
$\norm{\cdot}_{\widehat{K_f}}$.

Fixing $x \in \Real^n$, $x \neq 0$, we will show that $\norm{ \nabla
u(x) }_{\widehat{K_f}}^\ast  \leq C$ for some universal constant
$C>0$. Write $f = \exp(-g)$, and denote for short:
\[
A = \int_0^1 r^{n-1} f(rx) dr \text{ and } B = \int_1^{\infty}
r^{n-1} f(rx) dr~;
\]
note that:
\begin{multline*}
u = \left(A/(A+B)\right)^{1/n}~, \\
\nabla A = -\int_0^1 r^{n} f(rx) \nabla g(rx) dr~,
    \text{ and } \nabla B = -\int_1^{\infty} r^{n} f(rx) \nabla g(rx) dr~.
\end{multline*}
By Proposition \ref{prop:Fradelizi}, since $f(0)=1$ and $0$ is the
barycenter of $f$, then $\sup_x f(x) \leq \exp(n)$. This clearly
implies that $A \leq \exp(n)/n$, and that $g(x) \geq -n$. Denote
also:
\[
A^\ast = \int_0^1 r^{n} f(rx) \norm{\nabla
g(rx)}_{\widehat{K_f}}^\ast dr \text{ and } B^\ast = \int_1^{\infty}
r^{n} f(rx) \norm{\nabla g(rx)}_{\widehat{K_f}}^\ast dr.
\]
Then by (\ref{eq:u-definition})
\begin{equation}\begin{split}\label{eq:Lip-start}
\norm{ \nabla u(x) }_{\widehat{K_f}}^\ast &= \frac{1}{n}
\brac{\frac{A}{A+B}}^{\frac{1}{n}-1}
\frac{\norm{ \nabla A (A+B) - A (\nabla A + \nabla B) }_{\widehat{K_f}}^\ast}{(A+B)^2}  \\
&\leq \frac{1}{n} \brac{\frac{A}{A+B}}^{\frac{1}{n}}
\frac{A^\ast B + A B^\ast}{A (A+B)} \\
&\leq \frac{1}{n} \frac{A^\ast}{A} + \frac{1}{n}
\brac{\frac{A}{A+B}}^{\frac{1}{n}}
\frac{A^\ast+B^\ast}{A+B} \\
&\leq \frac{1}{n} \frac{A^\ast}{A} + \frac{e}{n}
\frac{A^\ast+B^\ast}{(A+B)^{\frac{n+1}{n}}}~.
\end{split}\end{equation}
Note that by the convexity of $g$, for all $x,y \in \Real^n$:
\[
g(y) \geq g(x) + \scalar{\nabla g(x) , y - x}.
\]
Recall the definition (\ref{eq:def-K_f^0}), stating that $y \in K_f^0$
iff $g(y) \leq n + g(0) = n$, and also recall that $g(x) \geq -n$.
This implies that for $y \in K_f^0$:
\[
\scalar{\nabla g(x) , y} \leq \scalar{\nabla g(x) , x} + g(y) - g(x)
\leq \scalar{\nabla g(x) , x} + 2n.
\]
By Proposition \ref{prop:K_f^0-relation} $\widehat{K_f} \subset K_f \subset
D K_f^0$, where $D = C (\sup_x f(x))^{1/n} \leq C e$ for some
universal $C>1$; hence
\[
\norm{\nabla g(x)}_{\widehat{K_f}}^\ast \leq D (\scalar{\nabla g(x)
, x} + 2n).
\]
We will use this rough estimate to bound $A^\ast$ and $B^\ast$ from
above. More generally, for $0\leq a<b \leq \infty$,
\begin{equation}
\begin{split}\label{eq:grad-bound}
\frac{1}{D} &\int_a^b r^{n} f(rx) \norm{\nabla
g(rx)}_{\widehat{K_f}}^\ast dr
    \leq \int_a^b r^{n} f(rx) (\scalar{\nabla g(rx) , rx} + 2n) dr \\
    &= \left . \frac{d}{dt} \right |_{t=1}
        \brac{-\int_a^b r^{n} f(t r x) dr} + 2n \int_a^b r^{n} f(rx) dr \\
    &= \left . \frac{d}{dt} \right |_{t=1} \brac{- t^{-(n+1)}
        \int_{at}^{bt} r^{n} f(r x) dr} + 2n \int_a^b r^{n} f(rx) dr \\
    &= (3n+1) \int_a^b r^{n} f(rx) dr + a^{n+1}
        f(ax) - b^{n+1} f(bx).
\end{split}\end{equation}
Of course the last term is interpreted as 0 when $b = \infty$. With
this bound in mind, let:
\[
A' = \int_0^1 r^{n} f(rx) dr \text{ and } B' = \int_1^\infty r^{n}
f(rx) dr~.
\]
Applying (\ref{eq:grad-bound}), we see that:
\begin{eqnarray*}
A^\ast / D &\leq& (3n+1) A' - f(x)~; \\
(A^\ast + B^\ast) / D &\leq& (3n+1) (A' + B')~,
\end{eqnarray*}
Hence by (\ref{eq:Lip-start})
\[
\norm{ \nabla u(x) }_{\widehat{K_f}}^\ast
    \leq \frac{(3n+1)D}{n}
        \brac{\frac{A'}{A} + e \frac{A'+B'}{(A+B)^{\frac{n+1}{n}}}}.
\]
Obviously $A' \leq A$ since $r \leq 1$ in the integrand of $A'$. By
Lemma \ref{lem:moments-lemma} (that is applicable since $f(0) = 1$)
we have:
\[
A'+B' = \int_0^\infty r^n f(xr) dr
    \leq C \brac{\int_0^\infty r^{n-1} f(xr) dr}^{\frac{n+1}{n}}
    = C (A+B)^{\frac{n+1}{n}},
\]
where $C>0$ is some universal constant. It follows that:
\[
\norm{u}_{Lip} = \sup_{x\in \Real^n} \norm{\nabla
u(x)}_{\widehat{K_f}}^\ast \leq 4D(1+e C).
\]

\medskip
\noindent{\bf Step 3}: Proof for general smooth functions.
\medskip

We have shown the assertion of the theorem for smooth functions $f$
with $f(0)=1$. In the general case, obviously $f(0)>0$, since the
barycenter of the log-concave $f$ is at the origin. Let us push
forward $f(x) dx$ by the map $S(x) = f(0)^{1/n} x$ to obtain $f'(x)
dx$, where:
\[
f'(x) = f(0)^{-1} f(f(0)^{-1/n} x) ~.
\]
Clearly $K_{f'}$ is a homothetic copy of $K_f$, and since
\[
\mes_n(K_{f'}) = \int f'(x) dx = \int f(x) dx = \mes_n(K_f)~,
\]
we see that $K_{f'} = K_f$. Let $T$ denote the radial map pushing
forward $f'(x) dx$ to the restriction of the Lebesgue measure on
$K_f$, denoted $\lambda$. Let $u' : (\Real^n,\norm{\cdot}_{K_f})
\rightarrow [0,1]$ be defined by:
\[
T'(x) = u'(x) \frac{x}{\norm{x}_{K_f}},
\]
and $u'(0)=0$. Since $f'(0)=1$ and $f'$ is smooth, step 2 implies
that $\norm{u'}_{Lip} \leq C$. Obviously $T = T' \circ S$ (e.g. by
uniqueness of the radial map pushing forward $f(x) dx$ onto
$\lambda$), and hence $u = u' \circ S$. This implies:
\[
\norm{u}_{Lip} = \norm{u'}_{Lip} f(0)^{1/n} \leq C f(0)^{1/n}~,
\]
and concludes the proof.
\qed

\begin{rem}
Of course the proof uses the fact that the barycenter of $f$ is at
the origin in a very indirect way. In fact, it is clear from the
proof that we may use any log-concave function $f$ for which:
\[
f(0) \geq D^{-n} \sup_{x \in \Real^n} f(x),
\]
for some $D \geq 1$, yielding 
$\norm{u}_{Lip} \leq C(D) f(0)^{1/n}$, where $C(D)$ is a constant
depending on $D$.
\end{rem}

As an immediate corollary of Theorem \ref{thm:Lip}, we obtain the
Bobkov-Ledoux Proposition from the introduction, although the direct
route taken by Bobkov and Ledoux in \cite{BobkovLedoux} is simpler
in this case and recovers a better universal constant in the bound.

\begin{proof}[Proof of the Bobkov--Ledoux Proposition]
It is easy to see that the Lipschitz constant of $S$ as a map acting
on $(\Real^n,\norm{\cdot})$ is invariant to scaling of the
Lebesgue measure, so we may assume that $\mes_n(K) = 1$. By
Theorem \ref{thm:Lip},
\[ \norm{S}_{Lip} \leq C f(0)^{1/n}
    = C \Gamma(1+n/p)^{-1/n}~.\]

\end{proof}

We will see in the next section how Theorem \ref{thm:Lip} may be
used to transfer isoperimetric inequalities from log-concave
measures to uniform measures on convex bodies.

\section{General Uniformly Convex Bodies}\label{sec:gucb}

In this section we give a proof of Proposition~\ref{prop:FigielPisier}
and provide the details that lead to Theorem~\ref{isop3}.

Let $\delta = \delta_V$ denote the modulus of convexity of a normed
space $V = (X,\norm{\cdot})$. It is known that $\delta$ is not
necessarily a convex function; we denote by $\tilde{\delta}$ the
maximal convex function majorated by $\delta$. We summarise several
known facts about $\delta$ and $\tilde{\delta}$ (see Lindenstrauss and
Tsafriri \cite[Proposition 1.e.6,Lemmata 1.e.7,1.e.8]{LT-Book-II}).

\begin{lem} \label{lem:delta-facts}
\hfill
\begin{enumerate}
\item
$\delta(t)/t$ is non-decreasing on $[0,2]$.
\item
$\delta(t/2) \leq \tilde{\delta}(t) \leq \delta(t)$ for all $t \in
[0,2]$.
\item
There exists a constant $C \geq 1$ such that
${\tilde{\delta}(t)}/{t^2} \leq C {\tilde{\delta}(s)}/{s^2}$, for
all $0 \leq t \leq s \leq 2$.
\end{enumerate}
\end{lem}

The following crucial fact is due to Figiel and Pisier \cite{FigielPisier}
(see also \cite[Lemma 1.e.10]{LT-Book-II}):

\begin{prop}[Figiel--Pisier] \label{lem:FigielPisier}
Let $x,y \in X$ such that $\norm{x}^2 + \norm{y}^2 = 2$. Then:
\[
\norm{x+y}^2 \leq 4 - 4 \delta(\norm{x-y}/2).
\]
\end{prop}

Proposition \ref{prop:FigielPisier} is an easy corollary of these
lemmata.

\begin{proof}[Proof of Proposition \ref{prop:FigielPisier}]
Let $x,y \in X$ such that $\norm{x}^2 + \norm{y}^2 \leq 2$, and
denote $s^2 := (\norm{x}^2 + \norm{y}^2)/2 \leq 1$. If $s = 0$ then
$\norm{x} = \norm{y} = 0$ and the claim is trivial. Otherwise,
denote $x' = x/s$ and $y' = y/s$, so that $\norm{x'}^2 + \norm{y'}^2
= 2$. Hence by Proposition \ref{lem:FigielPisier}:
\[
\norm{\frac{x'+y'}{2}}^2 \leq 1 -
\delta\brac{\frac{\norm{x'-y'}}{2}},
\]
or equivalently:
\[
\norm{\frac{x+y}{2}}^2 \leq s^2 - s^2
\delta\brac{\frac{\norm{x-y}}{2s}}.
\]
Now, $s \leq 1$; hence by Lemma \ref{lem:delta-facts} we have for
any $t \in [0,2s]$:
\[
s^2 \delta(t/s) \geq s^2 \tilde{\delta}(t/s) \geq c
\tilde{\delta}(t) \geq c \delta(t/2),
\]
where $c>0$ is a universal constant. Applying this for
\[
t = \frac{\norm{x-y}}{2} \leq \frac{\norm{x} + \norm{y}}{2} \leq s~,
\]
we conclude that:
\begin{equation} \label{eq:remark}
\norm{\frac{x+y}{2}}^2 \leq \frac{\norm{x}^2 + \norm{y}^2}{2}  - c
\delta\brac{\frac{\norm{x-y}}{4}},
\end{equation}
as required.
\end{proof}
\begin{rem}
Using $\norm{x}=\norm{y}=1$ in (\ref{eq:remark}), we see that
$\delta_V(\eps) \geq \frac{c}{2}\delta(\eps/4)$ for any function
$\delta$ satisfying (\ref{eq:remark}), so Proposition
\ref{prop:FigielPisier} is in fact a characterization (up to
universal constants) of the modulus of convexity $\delta_V$.
\end{rem}

Now we can fill the details in the proof of Theorem \ref{isop3}.
Assume that $V=(\Real^n,\norm{\cdot})$ is a uniformly convex space,
and let $\delta = \delta_V$ denote its modulus of convexity as
before. Scale the Lebesgue measure on $\Real^n$ so that
$\mes_n\set{\norm{x} \leq 1} = 1$, since the statement of Theorem
\ref{isop3} is invariant to this scaling. Now denote by $\mu$
the probability measure with density:
\[
f(x) = \frac{1}{Z} \exp(-n/c \norm{4 x}^2) \mathbf{1} (\norm{x}\leq 1/4)
\]
with respect to the Lebesgue measure, where $c>0$ is the constant
from Proposition \ref{prop:FigielPisier}. Here $Z>0$ is a
scaling factor so that $\mu$ be indeed a probability measure.
Integrating on level sets of $\norm{\cdot}$, it is clear that:
\[\begin{split}
Z &= \int_{\Real^n} \exp\brac{-\frac{1}{c} n \norm{4x}^2}
\mathbf{1} \brac{\norm{x}\leq \frac{1}{4}} dx \\
&= n \int_0^{1/4} \exp\brac{-\frac{16}{c} n s^2} s^{n-1} ds,
\end{split}\]
and in particular $Z^{1/n} \geq c' > 0$.

Write $f = \exp(-g)$, with
$g: \Real^n \rightarrow \Real \cup \set{+\infty}$.
Proposition \ref{prop:FigielPisier} then implies
that $g$ is uniformly convex, and satisfies:
\[
\frac{g(x) + g(y)}{2} - g\brac{\frac{x+y}{2}} \geq n \delta_1
(\norm{x-y}),
\]
where $\delta_1$ coincides with $\delta$ on $[0,1/4]$ and
$\delta_1(t) = +\infty$ for $t > 1/4$. Since $\delta(t)/t$
is non-decreasing by Lemma \ref{lem:delta-facts}, so is
$\delta_1(t)/t$, and assumption (\ref{eq:delta-assumption}) is
fulfilled. We can therefore apply Theorem \ref{isop1}, and deduce
an isoperimetric inequality for $\mu$ on $V$:
\[
\mu^+_{\| \cdot \|} (A)
    \geq C_{n,\delta} \widetilde{\mu(A)} \,
        \gamma_{n} \left(\log \frac{1}{\widetilde{\mu(A)}}\right)
        \text{ for all } A \subset \RR^n,
\]
where $C_{n,\delta}$ is given by $(\ref{eq:C-delta})$ and:
\[
\gamma_{n}(t) = \frac{t}{\delta_1^{-1}(t/(2 n))}~.
\]

We would now like to transfer this isoperimetric inequality to
$\lambda_V$, the uniform probability measure on $K_V = \set{\norm{x}\leq 1}$,
via a radial Lipschitz map. Clearly, $K_f$ is a homothetic copy of $K_V$, and since
\[ \mes_n(K_f) = \int f(x) dx = 1 = \mes_n(K_V)~,\]
it follows that $K_f = K_V$. Note also that
\[ f(0)^{1/n} = Z^{-1/n} \leq (c')^{-1}~.\]
Applying Theorem \ref{thm:Lip}, it follows that the Lipschitz
constant of the radial map pushing forward $\mu$ onto $\lambda_V$ is
bounded by a universal constant.
Because of the truncation in the definition of $\delta_1$, this only
implies the statement of Theorem \ref{isop3} for sets $A$ such that
\[ \widetilde{\lambda(A)} \geq \exp(-2 \delta(1/4) n)~. \]

Now suppose
\[ \widetilde{\lambda(A)} < \exp(-2 \delta(1/4) n)~. \]
Then
\[ \delta^{-1} \left( \frac{1}{2n} \log \frac{1}{\widetilde{\lambda(A)}} \right)
    \geq 1/4~,\]
and hence by Bobkov's inequality (\ref{bisop}) with $r = 1$
\[ \lambda_{\| \cdot \|}^+(A)
    \geq \frac{1}{2} \widetilde{\lambda}(A) \log \frac{1}{\widetilde{\lambda}(A)}
    \geq c' C_{n,\delta} \frac{\widetilde{\lambda}(A) \log \frac{1}{\widetilde{\lambda}(A)}}
                              {\delta^{-1} \left( \frac{1}{2n}
                                \log \frac{1}{\widetilde{\lambda(A)}} \right)} \]
with, say, $c' = e/(4(e-1))$~.

This concludes the proof of Theorem \ref{isop3}.


\section{Concentration and functional inequalities}\label{sec:5}

\subsection{Concentration of measure on uni\-form\-ly con\-vex bo\-dies}\label{subsec:conc}

In this subsection, we discuss the connection between our results and the
following Gromov--Milman inequality \cite{Gromov-Milman}, that we cite in the form
of Arias-de-Reyna, Ball, and Villa \cite{ABV}.

\begin{thm*}[Gromov--Milman]
Let $V = (\RR^n, \| \cdot \|)$ be a normed space; let $\delta = \delta_V$ be its
modulus of convexity, and let $\lambda$ be the uniform measure on the unit ball
of $V$. Then
\begin{equation}\label{eq:gm1}
1 - \lambda(B_{\eps,\| \cdot \|}) \leq \frac{1}{\lambda(B)} \exp(-2
n \delta(\eps)) \; \text{ for all } \; B \subset \RR^n~.
\end{equation}
In particular, if $\delta(\eps) \geq \alpha' \eps^p$, then
\begin{equation}\label{eq:gm2}
1 - \lambda(B_{\eps,\| \cdot \|}) \leq \frac{1}{\lambda(B)} \exp(-2
\alpha' n \eps^p)
    \; \text{ for all } \; B \subset \RR^n~.
\end{equation}
\end{thm*}

Let us compare this to our results. First assume $\delta(\eps) \geq
\alpha' \eps^p$; then (\ref{eq:puc_norm}) holds with $\alpha =
\alpha'/2^p$ (as mentioned in Subsection
\ref{sec:intro-uniformly-convex-bodies}). Therefore by
Theorem~\ref{isop2}
\begin{equation}\label{isop'}
\lambda_{\| \cdot \|}^+ (A) \geq C' (\alpha')^{1/p} n^{1/p}
    \widetilde{\lambda(A)} \log^{1-1/p} \frac{1}{\widetilde{\lambda(A)}}
    \;\; \text{ for all } \; B \subset \RR^n~,
\end{equation}
where $C'$ is a universal constant. Hence by Corollary \ref{cor:1.9}
\begin{equation}\label{eq:gm2'}
1 - \lambda(B_{\eps, \| \cdot \|})
    \leq \exp \left\{ -
        \left[ \log^{1/p} \frac{1}{1 - \lambda(B)} +
            \frac{c(\alpha')^{1/p}n^{1/p}\eps}{p}\right]^p
    \right\}
\end{equation}
The right-hand side in (\ref{eq:gm2'}) is at most
\[ (1 - \lambda(B)) \exp \left\{ - C' (\alpha')^{1/p} n^{1/p}
    \log^{1-1/p} \frac{1}{1 - \lambda(B)} \; \eps \right\}
    < 1 - \lambda(B)~;\]
hence (\ref{eq:gm2'}) yields a meaningful bound for any $\eps > 0$,
whereas (\ref{eq:gm2}) is meaningful for
\[ \eps \geq \left\{ \frac{1}{2\alpha'n}
     \log \frac{1}{\lambda(B)(1-\lambda(B))} \right\}^{1/p}~.\]

On the other hand, for larger $\eps$ the right-hand side of (\ref{eq:gm2'})
behaves like
\[ \exp\left\{ - \frac{C'^p}{p^p} \alpha' n \eps^p\right\}~;\]
that is, we lose a factor $p^p$ in the exponent.

The preceding discussion can be extended to arbitrary moduli of
convexity. In the general case, Theorem~\ref{isop3} yields
\begin{equation}\label{eq:gm1''}
\lambda_{\| \cdot \|}^+ (A) \geq
    C'_{n,\delta} \frac{\widetilde{\lambda(A)} \log \frac{1}{\widetilde{\lambda(A)}}}
                       {\delta^{-1}\left(\frac{1}{2n} \log \frac{1}{\widetilde{\lambda(A)}}
                            \right)}~;
\end{equation}
hence by Proposition~\ref{propA}
\begin{equation} \label{eq:propA-again}
1 - \lambda(B_{\eps, \| \cdot \|})
    \leq \exp \left\{ - h_{1 - \lambda(B)}^{-1}(\eps) \right\}~,
\end{equation}
where
\[ h_a(x) = \int_{\log 1/a}^x \frac{\delta^{-1}(y/2n) dy}{C_{n,\delta}' y}~. \]

By Lemma~\ref{lem:delta-facts} we can assume without loss of generality that
$\delta$ is convex (and $\delta^{-1}$ is concave). Then,
\[\begin{split}
h_a(x) &= \int_{\log 1/a}^x \frac{dy}{C_{n,\delta}' y} \,
    \frac{\int_{\log 1/a}^x \frac{\delta^{-1}(y/2n) dy}
         {C_{n,\delta}'y}}{\int_{\log 1/a}^x \frac{dy}{C_{n,\delta}' y}} \\
       &\leq \int_{\log 1/a}^x \frac{dy}{C_{n,\delta}' y} \,\,
            \delta^{-1} \left\{ \frac{1}{2n}
                \frac{\int_{\log 1/a}^x \frac{dy}{C_{n,\delta}'}}
                     {\int_{\log 1/a}^x \frac{dy}{C_{n,\delta}' y}} \right\} \\
       &= \frac{\log x - \log \log 1/a}{C_{n,\delta}'}
            \delta^{-1} \left\{ \frac{1}{2n}
                \frac{x - \log 1/a}{\log x - \log \log 1/a} \right\}~.
\end{split}\]
Now, $t \mapsto \delta^{-1}(t)/t$ is decreasing, hence
\begin{equation} \label{eq:GM-discuss}
h_a(x) \leq \frac{1}{C_{n,\delta}'}
    \delta^{-1} \left\{ \frac{1}{2n} (x - \log 1/a) \right\}
    \quad \text{if $x \leq e \log 1/a$.}
\end{equation}
On the other hand,
\[
h_a(e \log 1/a) = \int_{\log 1/a}^{e \log 1/a} \frac{\delta^{-1}(y/2n) dy}{C_{n,\delta}'y}
    \geq \frac{e-1}{eC_{n,\delta}'} \delta^{-1} \left\{ \frac{e \log 1/a}{2n}\right\}~;
\]
hence for $\eps \leq \frac{e-1}{eC_{n,\delta}'} \delta^{-1} \left\{
\frac{e \log 1/a}{2n}\right\}$, $x = h^{-1}_a(\eps) \leq e \log
1/a$, and (\ref{eq:GM-discuss}) implies:
\[
h_a^{-1}(\eps) \geq 2n\delta(C_{n,\delta}' \eps) + \log 1/a~.
\]
We conclude by (\ref{eq:propA-again}) that:
\begin{equation}\label{eq:gm1'}
1 - \lambda(B_{\eps, \| \cdot \|})
    \leq (1 - \lambda(B)) \exp \left\{ - 2n \delta(C_{n,\delta}' \eps)\right\}~.
\end{equation}

Again, (\ref{eq:gm1'}) is better than (\ref{eq:gm1}) for small $\eps$; if
\[ \eps \leq \frac{e-1}{eC_{n,\delta}'} \delta^{-1} \left\{ \frac{e \log 1/a}{2n}\right\}\]
the inequalities (\ref{eq:gm1}) and (\ref{eq:gm1'}) are similar, whereas for larger $\eps$
an inequality of type (\ref{eq:gm1'}) can only be deduced from (\ref{eq:gm1''}) under additional
regularity assumptions on $\delta$.

\subsection{Proofs}\label{subsec:concpr}

It remains to prove Propositions~\ref{propA} and \ref{propB}.

\begin{proof}[Proof of Proposition~\ref{propA}]
Let $B \subset \RR^n$ be a Borel set such that:
\[ a = 1 - \mu(B) \leq 1/2~;\]
the proof easily extends to the complementary case $a > 1/2$.

Denote $f(t) = 1 - \mu(B_t)$. Our assumptions then read:
\[ f(0) = a~; \quad df/dt(t) \leq -  f(t) \gamma(- \log f(t))\]
(where strictly speaking $df/dt$ should be the upper left
derivative). Setting $g = - \log f$,
\[ g(0) = \log 1/a~; \quad dg/dt \geq \gamma \circ g~,\]
and if $h = g^{-1}$,
\[ h(\log 1/a) = 0 \quad \text{and} \quad dh/dt \leq 1/(\gamma)~.\]
Therefore
\[ h(x) \leq \int_{\log 1/a}^x \frac{dy}{\gamma(y)} = h_a(x),\]
and
\[
f(t) = \exp(-h^{-1}(t)) \leq \exp(-h_a^{-1}(t)),
\]
as required.

The converse direction is obvious.
\end{proof}

\begin{proof}[Proof of Proposition~\ref{propB}]
Let us show that \ref{cap1}.\ implies \ref{cap2}. Let $F$ be a
function satisfying (\ref{techn}); assume for simplicity that the distribution
of $F$ has no atoms except for $0$ and $1$ and
that $\mu \{F = 0\} = 1/2$, $\mu \{ F = 1 \} = t = 1/2^k$. Choose
\[ 0 = u_1 < u_2 < \cdots < u_k = 1\]
so that
\[ \mu \{ u_i < F < u_{i+1} \} = 1/2^{i+1}~.\]
Then
\[\begin{split}
\int \| \nabla F \|_\ast^q d\mu
    &= \sum \int_{u_i < F \leq u_{i+1}} \| \nabla F \|_\ast^q d\mu \\
    &\geq \sum \frac{1}{2^{i+1}}
        \left\{ 2^{i+1} \int_{u_i < F < u_{i+1}}
            \| \nabla F \|_\ast d\mu \right\}^q
\end{split}\]
by Jensen's inequality. Now, apply \ref{cap1}.\ to the function
\[ F_i = \max\left( 0, \min \left( 1,
    \frac{F - u_i}{u_{i+1} - u_i} \right) \right)~.\]
Since $\mu\set{F_i = 1} = \mu\set{F \geq u_{i+1}} = 1/2^{i+1}$, we
obtain:
\[ \int_{u_i < F < u_{i+1}} \| \nabla F \|_\ast d\mu
    \geq c \, \mathbf{c}_0 (u_{i+1} - u_i) \frac{\log^{1/q} 2^{i+1}}{2^{i+1}}~; \]
therefore
\[\begin{split}
\int \| \nabla F \|_\ast^q d\mu
    &\geq \sum \frac{1}{2^{i+1}}
        \left\{ c \, \mathbf{c}_0 (u_{i+1} - u_i) \log^{1/q} 2^{i+1} \right\}^q \\
    &\geq c'' \, \mathbf{c}_0^q \sum (u_{i+1} - u_i)^q \, \frac{i+1}{2^{i+1}} \\
    &\geq c'' \, \mathbf{c}_0^q \left(\sum (u_{i+1} - u_i)\right)^q \Big/
        \left[ \sum \left( \frac{2^{i+1}}{i+1} \right)^{p/q} \right]^{q/p} \end{split}\]
according to H\"older's inequality. Finally,
\[\sum_{i=1}^k \left( \frac{2^{i+1}}{i+1} \right)^{p/q}
    \leq C (2^k / k)^{p/q} \]
and thence
\[\int \| \nabla F \|_\ast^q d\mu
     \geq c''' \, \mathbf{c}_0^q \, \frac{k}{2^k}
    \geq c' \, \mathbf{c}_0^q \, t \log 1/t~.\]
\end{proof}

\bibliographystyle{plain}
\bibliography{ConvexBib}

\def\cprime{$'$}
\begin{thebibliography}{10}

\bibitem{AleskerGM}
S.~Alesker.
\newblock Localization technique on the sphere and the {G}romov-{M}ilman
  theorem on the concentration phenomenon on uniformly convex sphere.
\newblock In {\em Convex geometric analysis (Berkeley, CA, 1996)}, volume~34 of
  {\em Math. Sci. Res. Inst. Publ.}, pages 17--27. Cambridge Univ. Press,
  Cambridge, 1999.

\bibitem{ABV}
J.~Arias-de Reyna, K.~Ball, and R.~Villa.
\newblock Concentration of the distance in finite-dimensional normed spaces.
\newblock {\em Mathematika}, 45(2):245--252, 1998.

\bibitem{BakryEmery}
D.~Bakry and M.~{\'E}mery.
\newblock Diffusions hypercontractives.
\newblock In {\em S\'eminaire de probabilit\'es, XIX, 1983/84}, volume 1123 of
  {\em Lecture Notes in Math.}, pages 177--206. Springer, Berlin, 1985.

\bibitem{BakryLedoux}
D.~Bakry and M.~Ledoux.
\newblock L\'evy-{G}romov's isoperimetric inequality for an
  infinite-dimensional diffusion generator.
\newblock {\em Invent. Math.}, 123(2):259--281, 1996.

\bibitem{Ball-PhD}
K.~Ball.
\newblock PhD thesis, Cambridge, 1986.

\bibitem{Ball-kdim-sections}
K.~Ball.
\newblock Logarithmically concave functions and sections of convex sets in
  $\mathbb{R}^n$.
\newblock {\em Studia Math.}, 88(1):69--84, 1988.

\bibitem{BallCarlenLieb}
K.~Ball, E.~A. Carlen, and E.~H. Lieb.
\newblock Sharp uniform convexity and smoothness inequalities for trace norms.
\newblock {\em Invent. Math.}, 115(3):463--482, 1994.

\bibitem{BartheGAFA}
F.~Barthe.
\newblock Log-concave and spherical models in isoperimetry.
\newblock {\em Geom. Funct. Anal.}, 12(1):32--55, 2002.

\bibitem{BartheRoberto}
F.~Barthe and C.~Roberto.
\newblock Sobolev inequalities for probability measures on the real line.
\newblock {\em Studia Math.}, 159(3):481--497, 2003.

\bibitem{BobkovExtremalHalfSpaces}
S.~Bobkov.
\newblock Extremal properties of half-spaces for log-concave distributions.
\newblock {\em Ann. Probab.}, 24(1):35--48, 1996.

\bibitem{BobkovGaussianIsoLogSobEquivalent}
S.~G. Bobkov.
\newblock Isoperimetric and analytic inequalities for log-concave probability
  measures.
\newblock {\em Ann. Probab.}, 27(4):1903--1921, 1999.

\bibitem{BobkovLocalizedProofOfGaussianIso}
S.~G. Bobkov.
\newblock A localized proof of the isoperimetric {B}akry-{L}edoux inequality
  and some applications.
\newblock {\em Teor. Veroyatnost. i Primenen.}, 47(2):340--346, 2002.

\bibitem{BobkovHoudre}
S.~G. Bobkov and C.~Houdr{\'e}.
\newblock Isoperimetric constants for product probability measures.
\newblock {\em Ann. Probab.}, 25(1):184--205, 1997.

\bibitem{BobkovLedoux}
S.~G. Bobkov and M.~Ledoux.
\newblock From {B}runn-{M}inkowski to {B}rascamp-{L}ieb and to logarithmic
  {S}obolev inequalities.
\newblock {\em Geom. Funct. Anal.}, 10(5):1028--1052, 2000.

\bibitem{BobkovZegarlinski}
S.~G. Bobkov and B.~Zegarlinski.
\newblock Entropy bounds and isoperimetry.
\newblock {\em Mem. Amer. Math. Soc.}, 176(829):x+69, 2005.

\bibitem{Borell-logconcave}
Ch. Borell.
\newblock Convex measures on locally convex spaces.
\newblock {\em Ark. Mat.}, 12:239--252, 1974.

\bibitem{Borell-GaussianIsoperimetry}
Ch. Borell.
\newblock The {B}runn--{M}inkowski inequality in {G}auss spaces.
\newblock {\em Inventiones Mathematicae}, 30:207--216, 1975.

\bibitem{Clarkson}
J.~A. Clarkson.
\newblock Uniformly convex spaces.
\newblock {\em Trans. Amer. Math. Soc.}, 40(3):396--414, 1936.

\bibitem{FedererFleming}
H.~Federer and W.~H. Fleming.
\newblock Normal and integral currents.
\newblock {\em Ann. of Math. (2)}, 72:458--520, 1960.

\bibitem{FigielModulusOfConvexity}
T.~Figiel.
\newblock An example of infinite dimensional reflexive {B}anach space
  non-isomorphic to its {C}artesian square.
\newblock {\em Studia Math.}, 42:295--306, 1972.

\bibitem{FigielPisier}
T.~Figiel and G.~Pisier.
\newblock S\'eries al\'eatoires dans les espaces uniform\'ement convexes ou
  uniform\'ement lisses.
\newblock {\em C. R. Acad. Sci. Paris S\'er. A}, 279:611--614, 1974.

\bibitem{FradeliziCentroid}
M.~Fradelizi.
\newblock Sections of convex bodies through their centroid.
\newblock {\em Arch. Math. (Basel)}, 69(6):515--522, 1997.

\bibitem{FradeliziGuedonLocalization1}
M.~Fradelizi and O.~Gu{\'e}don.
\newblock The extreme points of subsets of {$s$}-concave probabilities and a
  geometric localization theorem.
\newblock {\em Discrete Comput. Geom.}, 31(2):327--335, 2004.

\bibitem{FradeliziGuedonLocalization2}
M.~Fradelizi and O.~Gu{\'e}don.
\newblock A generalized localization theorem and geometric inequalities for
  convex bodies.
\newblock {\em Adv. Math.}, 204(2):509--529, 2006.

\bibitem{Gromov}
M.~Gromov.
\newblock {\em Metric structures for {R}iemannian and non-{R}iemannian spaces},
  volume 152 of {\em Progress in Mathematics}.
\newblock Birkh\"auser Boston Inc., Boston, MA, 1999.

\bibitem{Gromov-Milman}
M.~Gromov and V.~D. Milman.
\newblock Generalization of the spherical isoperimetric inequality to uniformly
  convex {B}anach spaces.
\newblock {\em Compositio Math.}, 62(3):263--282, 1987.

\bibitem{Hanner}
O.~Hanner.
\newblock On the uniform convexity of {$L\sp p$} and {$l\sp p$}.
\newblock {\em Ark. Mat.}, 3:239--244, 1956.

\bibitem{KLS}
R.~Kannan, L.~Lov{\'a}sz, and M.~Simonovits.
\newblock Isoperimetric problems for convex bodies and a localization lemma.
\newblock {\em Discrete Comput. Geom.}, 13(3-4):541--559, 1995.

\bibitem{KlartagPerturbationsWithBoundedLK}
B.~Klartag.
\newblock On convex perturbations with a bounded isotropic constant.
\newblock {\em Geom. and Funct. Anal.}, 16(6):1274--1290, 2006.

\bibitem{KlartagMilmanLogConcave}
B.~Klartag and V.~D. Milman.
\newblock Geometry of log-concave functions and measures.
\newblock {\em Geom. Dedicata}, 112:169--182, 2005.

\bibitem{LT-Book-II}
J.~Lindenstrauss and L.~Tzafriri.
\newblock {\em Classical {B}anach spaces. {II}}, volume~97 of {\em Ergebnisse
  der Mathematik und ihrer Grenzgebiete [Results in Mathematics and Related
  Areas]}.
\newblock Springer-Verlag, Berlin, 1979.
\newblock Function spaces.

\bibitem{LSLocalizationLemma}
L.~Lov{\'a}sz and M.~Simonovits.
\newblock Random walks in a convex body and an improved volume algorithm.
\newblock {\em Random Structures Algorithms}, 4(4):359--412, 1993.

\bibitem{Mazya}
V.~G. Maz{\cprime}ja~[Maz{\cprime}ya].
\newblock Classes of domains and imbedding theorems for function spaces.
\newblock {\em Soviet Math. Dokl.}, 1:882--885, 1960.

\bibitem{Milman-Pajor-LK}
V.~D. Milman and A.~Pajor.
\newblock Isotropic position and interia ellipsoids and zonoids of the unit
  ball of a normed $n$-dimensional space.
\newblock In {\em Geometric Aspects of Functional Analysis}, volume 1376 of
  {\em Lecture Notes in Mathematics}, pages 64--104. Springer-Verlag,
  1987-1988.

\bibitem{PisierUniformlyConvexThm}
G.~Pisier.
\newblock Martingales with values in uniformly convex spaces.
\newblock {\em Israel J. Math.}, 20(3-4):326--350, 1975.

\bibitem{SudakovTsirelson}
V.~N. Sudakov and B.~S. Cirel{\cprime}son~[Tsirelson].
\newblock Extremal properties of half-spaces for spherically invariant
  measures.
\newblock {\em Zap. Nau\v cn. Sem. Leningrad. Otdel. Mat. Inst. Steklov.
  (LOMI)}, 41:14--24, 165, 1974.
\newblock Problems in the theory of probability distributions, II.

\end{thebibliography}

\end{document}